\documentclass[11pt,a4paper]{article}

\usepackage{inputenc}
\usepackage{amsmath}
\usepackage{bm}
\usepackage{bbold}
\usepackage{amsthm}
\usepackage{enumerate}

\usepackage{tikz, tikz-3dplot}\tikzset{x=1cm,y=1cm,z=1cm}

\usepackage{pgfplots}\pgfplotsset{compat=1.14}
\usetikzlibrary{perspective}

\usepackage[hyphens]{url}

\usepackage{hyperref}
\usepackage{breakurl}

\title{Tropical optimization technique in bi-objective project scheduling under temporal constraints}

\author{N. Krivulin\thanks{Faculty of Mathematics and Mechanics, Saint Petersburg State University, 28 Universitetsky Ave., St.~Petersburg, 198504, Russia, 
nkk@math.spbu.ru.}
\thanks{This work was supported in part by the Russian Foundation for Basic Research (grant No. 20-010-00145).}
}

\date{}

\newtheorem{theorem}{Theorem}
\newtheorem{lemma}[theorem]{Lemma}
\newtheorem{corollary}[theorem]{Corollary}

\theoremstyle{definition}
\newtheorem{example}{Example}

\begin{document}

\maketitle

\begin{abstract}
We consider a project that consists of a set of activities performed in parallel under constraints on their start and finish times, including start-finish precedence relationships, release start times, release end times, and deadlines. The problems of interest are to decide on the optimal schedule of the activities to minimize both the maximum flow-time over all activities, and the project makespan. We formulate these problems as bi-objective optimization problems in the framework of tropical mathematics which investigates the theory and applications of algebraic systems with idempotent operations and has various applications in management science and operations research. Then, the use of methods and techniques of tropical optimization allows to derive complete Pareto-optimal solutions of the problems in a direct explicit form ready for further analysis and straightforward computation. We discuss the computational complexity of the solution and give illustrative examples.
\\

\\

\end{abstract}

\section{Introduction}

Tropical optimization deals with optimization problems that are formulated and solved in terms of tropical (idempotent) mathematics which concentrates on the theory and applications of algebraic systems with idempotent operations. Methods and techniques of tropical mathematics find application in operations research, management science and other fields, where tropical optimization allows to provide new efficient solutions to both known and novel optimization problems of practical importance.

Since the pioneering works by \cite{Pandit1961New,Cuninghamegreen1962Describing,Giffler1963Scheduling,Hoffman1963Onabstract,Vorobjev1963Extremal} on tropical mathematics in the early 1960s, optimization problems have served to motivate and illustrate the study. Further research in succeeding decades were often concerned with the analysis and solution of optimization problems as well. The results obtained in the area are presented in a number of contributed papers and books, among which are monographs by \cite{Baccelli1993Synchronization,Kolokoltsov1997Idempotent,Golan2003Semirings,Fiedler2006Linear,Heidergott2006Maxplus,Gondran2008Graphs,Gavalec2015Decision}.

Many tropical optimization problems consist in minimizing or maximizing functions defined in the tropical mathematics setting on vectors over idempotent semifields (semirings with idempotent addition and invertible multiplication). In some cases, but not infrequently, the problems can be solved analytically to provide a complete direct solution in an explicit form under rather general assumptions. For other problems only numerical techniques are known, available for a specific semifield, in the form of computational algorithms that offer a particular solution or indicate that no solution exists. A brief overview of the tropical optimization problems under consideration and their solutions can be found, e.g., in \cite{Krivulin2015Multidimensional}.

Applications of tropical mathematics include project scheduling problems, which appear in various settings in many publications from the early works by \cite{Pandit1961New,Cuninghamegreen1962Describing,Giffler1963Scheduling} to more recent papers \cite{Bouquard2006Application,Goto2009Robust,Singh2014Efficient} and books \cite{Baccelli1993Synchronization,Heidergott2006Maxplus,Gondran2008Graphs}. Models and methods of tropical optimization proved to be well suited for solving deterministic, pure temporal project scheduling problems, also referred to as time-dependent, time-constrained or resource-unconstrained problems (see, e.g., \cite{Demeulemeester2002Project,Neumann2003Project,Vanhoucke2012Project}). These problems involve temporal constraints (precedence relations, release times, due dates, deadlines) and temporal criteria (makespan, maximum deviation from due dates, maximum flow-time), and do not entail direct cost dependencies and resource requirements.

While general project scheduling problems are normally $N\!P$-hard, the temporal problems, due to the lack of cost and resource constraints imposed, can normally be formulated as linear programs or graph (network) optimization problems. As a result, they are numerically solved by appropriate computational procedures of polynomial computational complexity, such as the Karmarkar and Floyd-Warshall algorithms.

In contrast to the algorithmic solutions, the tropical optimization approach can offer complete, direct solutions to the problems, which are obtained in a compact vector form ready for further analysis and straightforward computations. Examples of temporal project scheduling problems and their direct analytical solutions in the framework of tropical optimization are provided, e.g., in \cite{Krivulin2015Extremal,Krivulin2017Direct,Krivulin2017Tropical,Krivulin2017Tropicaloptimization}.
 
Multi-objective project scheduling where two or more conflicting criteria have to be met \cite{Tkindt2006Multicriteria,Ballestin2015Handbook} allows the decision-maker to make more realistic and reasonable choice and thus is of particular importance. The common way to handle multi-objective problems is to obtain the best compromise solutions which constitute a set of non-dominated or Pareto-optimal solutions where no objective can be improved without degrading another one (see, e.g., \cite{Ehrgott2005Multicriteria,Luc2008Pareto,Miettinen2008Introduction,Benson2009Multiobjective}). For some bi-objective problems, a solution can be obtained analytically through the derivation of the Pareto frontier which is defined as the image of the Pareto-optimal set in the plane of objectives. The analytical description of the Pareto frontier is then used to find all Pareto optimal solutions \cite{Ruzika2005Approximation}. A natural approach to deal with multi-objective temporal project scheduling problems is to reduce them to multi-objective linear programs, which are then solved by appropriate numerical procedures like multiple objective variants of the simplex algorithm and the Benson algorithm. 

In this paper we consider a project that consists of a set of activities performed in parallel under constraints on their start and finish times, including start-finish precedence relationships, release start times, release end times, and deadlines. The problems of interest are to develop a schedule that minimizes both the the maximum flow-time over all activities and the project makespan, which present common objectives in project scheduling. We formulate and solve the problems in the framework of tropical mathematics as a tropical bi-objective optimization problem.

To handle the problem, we follow the approach developed in \cite{Krivulin2015Extremal,Krivulin2015Multidimensional,Krivulin2017Direct} to solve ordinary single objective problems and then applied to bi-objective problems in \cite{Krivulin2020Using}. The approach offers tropical optimization techniques that involve the introduction of parameters to represent the optimal values of objective functions and hence to reduce the optimization problem to a system of parametrized vector inequalities. We exploit the existence conditions for solutions of the system to evaluate the parameters and then to describe the Pareto frontier of the problem. The complete solution of the system, which corresponds to the parameters given by the Pareto frontier is taken as a Pareto-optimal solution of the initial bi-objective problem. We apply this result to solve the project scheduling problems of interest in a direct explicit form suitable for both formal analysis and numerical implementation. We discuss the computational complexity of the solution and give illustrative examples.

The paper is organized as follows. In Section~\ref{S-BOPSP} we provide a formal description of the bi-objective temporal project scheduling problems of interest. Section~\ref{S-PDNR} includes a brief overview of basic facts about tropical algebra, which are used in the subsequent solution of a bi-objective tropical optimization problem. The main result which provides a complete Pareto-optimal solution to the optimization problem in an exact analytical form is given in Section~\ref{S-BOTOP}. We apply the result obtained to solve the project scheduling problems and present illustrative examples in Section~\ref{S-ABOSP}.

\section{Bi-objective project scheduling problems}
\label{S-BOPSP}

We start with a formal description of bi-objective temporal project scheduling problems that serve to motivate and illustrate the solution obtained below in the framework of tropical optimization. In order to facilitate formulation of the problems in terms of tropical algebra, we use rather general model and notation which are slightly different from those commonly adopted in the literature on project scheduling (see, e.g., \cite{Demeulemeester2002Project,Neumann2003Project,Tkindt2006Multicriteria,Vanhoucke2012Project}).

Consider a project which consists of $n$ activities (jobs, tasks, operations) to be performed in parallel under certain temporal constraints. The problem is to construct a schedule for the activities, which can be considered optimal in the sense of the best compromise solutions with respect to two different criteria to be described later.

For each activity $i=1,\ldots,n$, we denote the unknown start time by $x_{i}$ and finish time by $y_{i}$, and assume these variables to be subject to constraints imposed due to some technological, organizational or other limitations of the project. First, we suppose that the start time may be bounded by a range given by specified release start time (ready or arrival time) $g_{i}$ and release deadline (release end time or time limit) $h_{i}$, which yields the double inequality constraints 
\begin{equation*}
g_{i}
\leq
x_{i}
\leq
h_{i},
\qquad
i=1,\ldots,n.
\end{equation*}

Furthermore, the minimum allowed time lag $a_{ij}$ between the start of activity $j$ and the finish of $i$ is given to specify the start-finish constraints by the inequalities
\begin{equation*}
a_{ij}+x_{j}
\leq
y_{i},
\qquad
i,j=1,\ldots,n,
\end{equation*}
where we take $a_{ij}=-\infty$ if the lag $a_{ij}$ is not defined. Note that the difference between the finish and start times of activity $i$ is bounded from below by the value of $a_{ii}$, which is normally assumed to be non-negative and presents the duration of the activity when no other activities are taken into account.

Moreover, we assume that each activity finishes immediately as soon as all its related start-finish constraints are satisfied, which results in the equalities
\begin{equation*}
\max_{1\leq j\leq n}(a_{ij}+x_{j})
=
y_{i},
\qquad
i=1,\ldots,n.
\end{equation*}

Finally, the finish time of activity $i$ may be bounded from above by a given deadline $f_{i}$ to satisfy the constraints
\begin{equation*}
y_{i}
\leq
f_{i},
\qquad
i=1,\ldots,n.
\end{equation*}

To develop an optimal schedule, we consider two optimality criteria which are frequently used in project scheduling and involve the minimization of the maximum flow-time of activities and the minimization of the project makespan. The flow-time (shop or cycle time) is defined for activity $i$ as the difference $y_{i}-x_{i}$ between its finish and start times, and may directly or indirectly reflect the expenditure incurred to perform the activity. The maximum flow-time over all activities is given by  
\begin{equation*}
\max_{1\leq i\leq n}(y_{i}-x_{i}).
\end{equation*}

The project makespan which is the overall project duration presents a commonly used measure of schedule efficiency to be minimized. This measure is calculated as the difference between the maximum finish time and the minimum start time of activities, it is not less than the maximum flow-time and represented as
\begin{equation*}
\max_{1\leq i\leq n}y_{i}-\min_{1\leq i\leq n}x_{i}
=
\max_{1\leq i\leq n}y_{i}+\max_{1\leq i\leq n}(-x_{i}).
\end{equation*}

The problems of interest are formulated to find the start and finish times for all activities to minimize both the maximum flow time and the makespan under the release time, release deadline, start-finish and deadline constraints. By combining the objective functions with the release and start-finish constraints, we arrive at the bi-objective problem of project scheduling in which, given numbers $a_{ij}$, $g_{i}$ and $h_{i}$ such that $g_{i}\leq h_{i}$ for all $i,j=1,\ldots,n$, one needs to obtain the unknown $x_{i}$ and $y_{i}$ to
\begin{equation}
\begin{aligned}
&
\text{minimize}
&&
\left\{
\max_{1\leq i\leq n}(y_{i}-x_{i}),\ 
\max_{1\leq i\leq n}y_{i}+\max_{1\leq i\leq n}(-x_{i})
\right\};
\\
&
\text{subject to}
&&
\max_{1\leq j\leq n}(a_{ij}+x_{j})
=
y_{i},
\\
&&&
g_{i}
\leq
x_{i}
\leq
h_{i},
\qquad
i=1,\ldots,n.
\label{P-minyixi-yixi-aijxjyi-gixihi}
\end{aligned}
\end{equation}

Another problem which takes into account the deadline constraints $f_{i}$ instead of the release deadlines $h_{i}$ is formulated as follows:
\begin{equation}
\begin{aligned}
&
\text{minimize}
&&
\left\{
\max_{1\leq i\leq n}(y_{i}-x_{i}),\ 
\max_{1\leq i\leq n}y_{i}+\max_{1\leq i\leq n}(-x_{i})
\right\};
\\
&
\text{subject to}
&&
\max_{1\leq j\leq n}(a_{ij}+x_{j})
=
y_{i},
\\
&&&
g_{i}
\leq
x_{i},
\qquad
y_{i}
\leq
f_{i},
\qquad
i=1,\ldots,n.
\label{P-minyixi-yixi-aijxjyi-gixi-yifi}
\end{aligned}
\end{equation}

Below, we apply methods and techniques of tropical optimization to represent bi-criteria optimization problems \eqref{P-minyixi-yixi-aijxjyi-gixihi} and \eqref{P-minyixi-yixi-aijxjyi-gixi-yifi} in compact vector form, and then to derive complete direct Pareto-optimal solutions of the problems.

\section{Preliminary definitions, notation and results}
\label{S-PDNR}

In this section, we offer an overview of basic definitions, notation and preliminary results of tropical (idempotent) algebra to provide an appropriate analytical framework for compact formulation and complete solution of a tropical optimization problem in the next section. Further details on the theory and applications of tropical mathematics can be found, e.g., in the monographs and textbooks by \cite{Baccelli1993Synchronization,Kolokoltsov1997Idempotent,Golan2003Semirings,Heidergott2006Maxplus,Gondran2008Graphs}.

\subsection{Idempotent semifield}

Let $\mathbb{X}$ be a set closed under two associative and commutative operations, addition $\oplus$ and multiplication $\otimes$, which have neutral elements, zero $\mathbb{0}$ and unit $\mathbb{1}$. Addition is idempotent, that is, $x\oplus x=x$ for any $x\in\mathbb{X}$. Multiplication distributes over addition, has $\mathbb{0}$ as absorbing element and is invertible to endow each $x\ne\mathbb{0}$ with the inverse $x^{-1}$ such that $x\otimes x^{-1}=\mathbb{1}$. The system $(\mathbb{X},\mathbb{0},\mathbb{1},\oplus,\otimes)$ is normally referred to as the idempotent semifield.

Integer powers are introduced in the standard way to represent iterated products as $\mathbb{0}^{n}=\mathbb{0}$, $x^{n}=x\otimes x^{n-1}$ and $x^{-n}=(x^{-1})^{n}$ for any $x\ne\mathbb{0}$ and natural $n$. Moreover, the semifield is assumed algebraically complete in the sense that the equation $x^{n}=a$ is uniquely solvable in $x$ for any $a\in\mathbb{X}$ and natural $n$, which makes rational exponents well defined. The multiplication symbol $\otimes$ is henceforth omitted to save writing.

Idempotent addition defines a partial order on $\mathbb{X}$ by the rule that $x\leq y$ if and only if $x\oplus y=y$. In terms of this order, addition and multiplication are monotone in both operands, that is, the inequality $x\leq y$ results in the inequalities $x\oplus z\leq y\oplus z$ and $xz\leq yz$ for $x,y,z\in\mathbb{X}$. Inversion is antitone, which means that for all $x,y\ne\mathbb{0}$ the inequality $x\leq y$ yields $x^{-1}\geq y^{-1}$. Addition has an extremal property (the majority law) that the inequalities $x\leq x\oplus y$ and $y\leq x\oplus y$ hold for any $x,y\in\mathbb{X}$. Moreover, the inequality $x\oplus y\leq z$ is equivalent to the system of two inequalities $x\leq z$ and $y\leq z$. In what follows, the partial order is assumed extended to a consistent total order to make the semifield linearly ordered.  

For any $x,y\in\mathbb{X}$ and non-negative integer $m$, the idempotent analogue of binomial identity is given by $(x\oplus y)^{m}=x^{m}\oplus x^{m-1}y\oplus\cdots\oplus y^{m}$. 

It follows from the identity with $m=2$ that the inequality $(xy)^{1/2}\leq x\oplus y$ is valid as an idempotent analogue of the relation between geometric and arithmetic means of two positive numbers. This inequality readily extends to any integer $m>0$ and $x_{1},\ldots,x_{m}\in\mathbb{X}$ in the form $(x_{1}\cdots x_{m})^{1/m}\leq x_{1}\oplus\cdots\oplus x_{m}$.

An example of the idempotent semifield under consideration is the real semifield $\mathbb{R}_{\max,+}=(\mathbb{R}\cup\{-\infty\},-\infty,0,\max,+)$ which is frequently called the max-plus algebra. In the semifield, the addition $\oplus$ is computed as calculation of maximum and the multiplication $\otimes$ as arithmetic addition; the zero $\mathbb{0}$ is defined as $-\infty$ and the unit $\mathbb{1}$ as the arithmetic zero $0$. Furthermore, for any $x\in\mathbb{R}$, there exists the inverse $x^{-1}$ which corresponds to the opposite number $-x$ in conventional algebra. The power $x^{y}$ coincides with the arithmetic product $yx$ which is defined for all $x,y\in\mathbb{R}$. The order induced by idempotent addition conforms with the natural linear order on $\mathbb{R}\cup\{-\infty\}$.

As another example, consider $\mathbb{R}_{\min,\times}=(\mathbb{R}_{+}\cup\{\infty\},\infty,1,\min,\times)$ where $\mathbb{R}_{+}$ is the set of positive reals. In the semifield, the operations are defined as $\oplus=\min$ and $\otimes=\times$, and the neutral elements as $\mathbb{0}=\infty$ and $\mathbb{1}=1$. The multiplicative inverses and exponents are interpreted as usual, whereas the order produced by addition is opposite to the natural linear order on $\mathbb{R}_{+}$.

\subsection{Algebra of matrices and vectors}

Denote the set of matrices over $\mathbb{X}$ with $m$ rows and $n$ columns by $\mathbb{X}^{m\times n}$. A matrix that has all entries equal to $\mathbb{0}$ is the zero matrix. If a matrix has no zero columns, it is called column-regular. 

Addition and multiplication of matrices as well as multiplication of matrices by scalars follow the standard rules where the arithmetic addition and multiplication are replaced by $\oplus$ and $\otimes$.

The transpose of a matrix $\bm{A}$ is denoted $\bm{A}^{T}$. The multiplicative conjugate transpose of a nonzero $(m\times n)$-matrix $\bm{A}=(a_{ij})$ is the $(n\times m)$-matrix $\bm{A}^{-}=(a_{ij}^{-})$ with the entries $a_{ij}^{-}=a_{ji}^{-1}$ if $a_{ji}\ne\mathbb{0}$, and $a_{ij}^{-}=\mathbb{0}$ otherwise.

The properties of scalar addition and multiplication, which are associated with the order relation induced by idempotent addition extend to the matrix operations where the inequalities are considered entry-wise.

A square matrix that has all diagonal entries equal to $\mathbb{1}$ and all off-diagonal entries to $\mathbb{0}$ is the identity matrix denoted by $\bm{I}$. Natural powers of a square matrix $\bm{A}$ are given by $\bm{A}^{0}=\bm{I}$ and $\bm{A}^{n}=\bm{A}^{n-1}\bm{A}$ for any natural $n$.

Let $\bm{A}=(a_{ij})$ be a square matrix of order $n$. The trace of the matrix $\bm{A}$ is calculated as $\mathop\mathrm{tr}\bm{A}=a_{11}\oplus\cdots\oplus a_{nn}$. For any matrices $\bm{A}$ and $\bm{B}$ of appropriate dimensions, and scalar $x$, the following equalities hold:
\begin{equation*}
\mathop\mathrm{tr}(\bm{A}\oplus\bm{B})
=
\mathop\mathrm{tr}\bm{A}
\oplus
\mathop\mathrm{tr}\bm{B},
\qquad
\mathop\mathrm{tr}(\bm{A}\bm{B})
=
\mathop\mathrm{tr}(\bm{B}\bm{A}),
\qquad
\mathop\mathrm{tr}(x\bm{A})
=
x\mathop\mathrm{tr}\bm{A}.
\end{equation*}

Consider a function that maps any matrix $\bm{A}\in\mathbb{X}^{n\times n}$ to the scalar
\begin{equation*}
\mathop\mathrm{Tr}(\bm{A})
=
\mathop\mathrm{tr}\bm{A}
\oplus
\cdots
\oplus
\mathop\mathrm{tr}\bm{A}^{n}.%,
\end{equation*}

Provided that $\mathop\mathrm{Tr}(\bm{A})\leq\mathbb{1}$, one can define the matrix (also known as the Kleene star matrix) given by
\begin{equation*}
\bm{A}^{\ast}
=
\bm{I}
\oplus
\bm{A}
\oplus\cdots\oplus
\bm{A}^{n-1}.
\end{equation*}

Any matrix that consists of a single column (row) forms a column (row) vector. In the following, all vectors are taken as column vectors unless otherwise indicated. The set of column vectors of size $n$ is denoted $\mathbb{X}^{n}$. 

A vector without zero components is called regular. The zero vector and the vector with all components equal to $\mathbb{1}$ are denoted by $\bm{0}=(\mathbb{0},\ldots,\mathbb{0})^{T}$ and $\bm{1}=(\mathbb{1},\ldots,\mathbb{1})^{T}$.

The conjugate transpose of a nonzero column vector $\bm{x}=(x_{i})$ is the row vector $\bm{x}^{-}=(x_{i}^{-})$ where $x_{i}^{-}=x_{i}^{-1}$ if $x_{i}\ne\mathbb{0}$, and $x_{i}^{-}=\mathbb{0}$ otherwise.

A vector $\bm{y}\in\mathbb{X}^{n}$ is linearly dependent on $\bm{x}_{1},\ldots,\bm{x}_{m}\in\mathbb{X}^{n}$ if there exist scalars $c_{1},\ldots,c_{m}\in\mathbb{X}$ such that $\bm{y}=c_{1}\bm{x}_{1}\oplus\cdots\oplus c_{m}\bm{x}_{m}$. A vector $\bm{y}$ is collinear with $\bm{x}$ if $\bm{y}=c\bm{x}$ for some $c$. 

A scalar $\lambda$ is an eigenvalue of a square matrix $\bm{A}\in\mathbb{X}^{n\times n}$ if there exists a nonzero vector $\bm{x}\in\mathbb{X}^{n}$ that satisfies the equation
$\bm{A}\bm{x}=\lambda\bm{x}$. The maximum eigenvalue is referred to as the spectral radius and calculated as
\begin{equation*}
\lambda
=
\mathop\mathrm{tr}\bm{A}
\oplus
\cdots
\oplus
\mathop\mathrm{tr}\nolimits^{1/n}(\bm{A}^{n}).
\end{equation*}

For any matrix $\bm{A}\in\mathbb{X}^{n\times n}$ and vector $\bm{x}\in\mathbb{X}^{n}$, idempotent analogues of matrix and vector norms are given by
\begin{equation*}
\|\bm{A}\|
=
a_{11}\oplus\cdots\oplus a_{nn}
=
\bm{1}^{T}\bm{A}\bm{1},
\qquad
\|\bm{x}\|
=
x_{1}\oplus\cdots\oplus x_{n}
=
\bm{1}^{T}\bm{x}
=
\bm{x}^{T}\bm{1}.
\end{equation*}

\subsection{Vector inequalities}

Suppose that, for given matrix $\bm{A}\in\mathbb{X}^{m\times n}$ and vector $\bm{d}\in\mathbb{X}^{m}$, one needs to solve, with respect to the unknown vector $\bm{x}\in\mathbb{X}^{n}$, the inequality
\begin{equation}
\bm{A}\bm{x}
\leq
\bm{d}.
\label{I-Axleqd}
\end{equation}

A solution to inequality \eqref{I-Axleqd} is described as follows (see, e.g., \cite{Baccelli1993Synchronization,Heidergott2006Maxplus,Krivulin2015Extremal}).
\begin{lemma}
\label{L-Axleqd}
For any column-regular matrix $\bm{A}$ and regular vector $\bm{d}$, all solutions to inequality \eqref{I-Axleqd} are given by the inequality $\bm{x}\leq(\bm{d}^{-}\bm{A})^{-}$.
\end{lemma}

Now assume that, given a matrix $\bm{A}\in\mathbb{X}^{n\times n}$ and vector $\bm{b}\in\mathbb{X}^{n}$, the problem is to find regular vectors $\bm{x}\in\mathbb{X}^{n}$ to satisfy the inequality
\begin{equation}
\bm{A}\bm{x}\oplus\bm{b}
\leq
\bm{x}.
\label{I-Axbleqx}
\end{equation}

The next result obtained in \cite{Krivulin2015Multidimensional} offers a direct solution.
\begin{theorem}
\label{T-Axbleqx}
For any matrix $\bm{A}$, the following statements hold:
\begin{enumerate}
\item If $\mathop\mathrm{Tr}(\bm{A})\leq\mathbb{1}$, then all regular solutions to \eqref{I-Axbleqx} are given by $\bm{x}=\bm{A}^{\ast}\bm{u}$ where $\bm{u}\geq\bm{b}$.
\item If $\mathop\mathrm{Tr}(\bm{A})>\mathbb{1}$, then there is only the trivial solution $\bm{x}=\bm{0}$.
\end{enumerate}
\end{theorem}

\subsection{Binomial identities for matrices and traces}

We conclude the overview of preliminary results with some useful matrix formulas. We start with a binomial identity that is valid for any matrices $\bm{A},\bm{B}\in\mathbb{X}^{n\times n}$ and integer $m\geq0$ in the following form (see also \cite{Krivulin2017Direct}):
\begin{equation*}
(\bm{A}\oplus\bm{B})^{m}
=
\bm{A}^{m}
\oplus
\bigoplus_{k=1}^{m}\bigoplus_{\substack{i_{0}+i_{1}+\cdots+i_{k}=m-k\\i_{0},i_{1},\ldots,i_{k}\geq0}}
\bm{A}^{i_{0}}(\bm{B}\bm{A}^{i_{1}}\cdots\bm{B}\bm{A}^{i_{k}}).
\end{equation*}

In the case when $\bm{B}=\bm{p}\bm{q}^{-}$ where $\bm{p},\bm{q}\in\mathbb{X}^{n}$ are vectors, this identity becomes
\begin{equation}
(\bm{A}\oplus\bm{p}\bm{q}^{-})^{m}
\\=
\bm{A}^{m}
\oplus
\bigoplus_{k=1}^{m}\bigoplus_{\substack{i_{0}+i_{1}+\cdots+i_{k}=m-k\\i_{0},i_{1},\ldots,i_{k}\geq0}}
(\bm{q}^{-}\bm{A}^{i_{1}}\bm{p}\cdots\bm{q}^{-}\bm{A}^{i_{k-1}}\bm{p})
\bm{A}^{i_{0}}\bm{p}\bm{q}^{-}\bm{A}^{i_{k}},
\label{E-Apqm}
\end{equation}
where we assume $\bm{q}^{-}\bm{A}^{i_{1}}\bm{p}\cdots\bm{q}^{-}\bm{A}^{i_{k-1}}\bm{p}=\mathbb{1}$ if $k=1$.

After taking trace of \eqref{E-Apqm} and summing up the traces over all $m=1,\ldots,n$, we apply properties of trace to rearrange the terms and obtain
\begin{equation}
\mathop\mathrm{Tr}(\bm{A}\oplus\bm{p}\bm{q}^{-})
\\=
\bigoplus_{k=1}^{n}\mathop\mathrm{tr}\bm{A}^{k}
\oplus
\bigoplus_{k=1}^{n}
\bigoplus_{m=0}^{n-k}
\bigoplus_{\substack{i_{1}+\cdots+i_{k}=m\\i_{1},\ldots,i_{k}\geq0}}
(\bm{q}^{-}\bm{A}^{i_{1}}\bm{p}\cdots\bm{q}^{-}\bm{A}^{i_{k}}\bm{p}).
\label{E-TrAB}
\end{equation} 

Suppose that $\mathop\mathrm{Tr}(\bm{A}\oplus\bm{p}\bm{q}^{-})\leq\mathbb{1}$. Then, the evaluation of the sum of \eqref{E-Apqm} over all $m=0,\ldots,n-1$ and rearrangement of terms yield the identity for the Kleene matrix
\begin{equation}
(\bm{A}\oplus\bm{p}\bm{q}^{-})^{\ast}
=
\bm{A}^{\ast}
\oplus
\bigoplus_{k=1}^{n-1}
\bigoplus_{m=0}^{n-k-1}
\bigoplus_{\substack{i_{0}+i_{1}+\cdots+i_{k}=m\\i_{0},i_{1},\ldots,i_{k}\geq0}}
(\bm{q}^{-}\bm{A}^{i_{1}}\bm{p}\cdots\bm{q}^{-}\bm{A}^{i_{k-1}}\bm{p})
\bm{A}^{i_{0}}\bm{p}\bm{q}^{-}\bm{A}^{i_{k}}.
\label{E-Apqast}
\end{equation}

\section{Bi-objective tropical optimization problem}
\label{S-BOTOP}

We are now in a position to formulate and solve a new bi-objective tropical optimization problem which is used in the next section to solve the project scheduling problems of interest. For the sake of compactness and generality of presentation, we describe here the tropical optimization problem and its solution in terms of a general idempotent semifield $(\mathbb{X},\mathbb{0},\mathbb{1},\oplus,\otimes)$.  

Suppose that, given a nonzero square matrix $\bm{A}\in\mathbb{X}^{n\times n}$ and nonzero vectors $\bm{p},\bm{q},\bm{g},\bm{h}\in\mathbb{X}^{n}$ where $\bm{g}\leq\bm{h}$, we need to find regular vectors $\bm{x}\in\mathbb{X}^{n}$ that provide solutions to the bi-objective optimization problem
\begin{equation}
\begin{aligned}
&
\text{minimize}
&&
\{\bm{x}^{-}\bm{A}\bm{x},\ 
\bm{x}^{-}\bm{p}\bm{q}^{-}\bm{x}\};
\\
&
\text{subject to}
&&
\bm{g}
\leq
\bm{x}
\leq
\bm{h}.
\label{P-minxAxxpqx-gxh}
\end{aligned}
\end{equation}

To solve this tropical optimization problem, we apply the solution technique developed in \cite{Krivulin2020Using} for bi-criteria decision-making problems. We reduce problem \eqref{P-minxAxxpqx-gxh} to a system of parameterized vector inequalities, where the parameters are used to describe the Pareto frontier of the problem. The existence conditions for solutions of the system serve to derive the Pareto frontier, whereas all corresponding solutions of the system give the Pareto-optimal solutions of \eqref{P-minxAxxpqx-gxh}.

\subsection{Parametrization of problem}

We start with a description of the Pareto frontier on the border of the image of the set of feasible solutions in the plane of criteria. We denote the minimum values of the scalar objective functions $\bm{x}^{-}\bm{A}\bm{x}$ and $\bm{x}^{-}\bm{p}\bm{q}^{-}\bm{x}$ in the Pareto frontier of problem \eqref{P-minxAxxpqx-gxh} by $\alpha$ and $\beta$, and note that $\alpha,\beta>\mathbb{0}$ for any regular $\bm{x}$. Then, all solutions are defined by the system of inequalities
\begin{equation*}
\bm{x}^{-}\bm{A}\bm{x}
\leq
\alpha,
\qquad
\bm{x}^{-}\bm{p}\bm{q}^{-}\bm{x}
\leq
\beta,
\qquad
\bm{g}
\leq
\bm{x}
\leq
\bm{h}.
\end{equation*}

By using Lemma~\ref{L-Axleqd}, we solve the first inequality with respect to $\bm{A}\bm{x}$ and the second to $\bm{p}\bm{q}^{-}\bm{x}$ to obtain the equivalent system
\begin{equation*}
\alpha^{-1}\bm{A}\bm{x}
\leq
\bm{x},
\qquad
\beta^{-1}\bm{p}\bm{q}^{-}\bm{x}
\leq
\bm{x},
\qquad
\bm{g}
\leq
\bm{x}
\leq
\bm{h},
\end{equation*}
which can then be combined into one double inequality
\begin{equation}
(\alpha^{-1}\bm{A}
\oplus
\beta^{-1}\bm{p}\bm{q}^{-})
\bm{x}
\oplus
\bm{g}
\leq
\bm{x}
\leq
\bm{h}.
\label{I-alpha1Abeta1pqxgleqxleqh}
\end{equation}

According to Theorem~\ref{T-Axbleqx}, regular solutions of the left inequality at \eqref{I-alpha1Abeta1pqxgleqxleqh} exist if and only if the condition $\mathop\mathrm{Tr}(\alpha^{-1}\bm{A}
\oplus\beta^{-1}\bm{p}\bm{q}^{-})\leq\mathbb{1}$ is valid, under which all solutions are given through a vector of parameters $\bm{u}$ by
\begin{equation*}
\bm{x}
=
(\alpha^{-1}\bm{A}
\oplus
\beta^{-1}\bm{p}\bm{q}^{-})^{\ast}
\bm{u},
\qquad
\bm{u}
\geq
\bm{g}.
%\label{E-x-alpha1Abeta111A-ug}
\end{equation*}

Furthermore, to provide the right inequality at \eqref{I-alpha1Abeta1pqxgleqxleqh}, the vector $\bm{u}$ must satisfy the inequality $(\alpha^{-1}\bm{A}\oplus\beta^{-1}\bm{p}\bm{q}^{-})^{\ast}\bm{u}\leq\bm{h}$. We apply Lemma~\ref{L-Axleqd} to this inequality and obtain an upper bound on $\bm{u}$, which together with the lower bound $\bm{g}$ gives the inequality
\begin{equation*}
\bm{g}
\leq
\bm{u}
\leq
(\bm{h}^{-}
(\alpha^{-1}\bm{A}
\oplus
\beta^{-1}\bm{p}\bm{q}^{-})^{\ast})^{-}.
\end{equation*}

The set of vectors $\bm{u}$ defined by this double inequality is not empty if the inequality $\bm{g}\leq(\bm{h}^{-}(\alpha^{-1}\bm{A}\oplus\beta^{-1}\bm{p}\bm{q}^{-})^{\ast})^{-}$ holds, which is, due to Lemma~\ref{L-Axleqd}, equivalent to the condition $\bm{h}^{-}(\alpha^{-1}\bm{A}\oplus\beta^{-1}\bm{p}\bm{q}^{-})^{\ast}\bm{g}\leq\mathbb{1}$.

By collecting the existence conditions and combining the bounds on $\bm{u}$, we finally conclude that inequality \eqref{I-alpha1Abeta1pqxgleqxleqh} has regular solutions if and only if both inequalities
\begin{equation}
\mathop\mathrm{Tr}(\alpha^{-1}\bm{A}
\oplus
\beta^{-1}\bm{p}\bm{q}^{-})
\leq
\mathbb{1},
\qquad
(\bm{h}^{-}
(\alpha^{-1}\bm{A}
\oplus
\beta^{-1}\bm{p}\bm{q}^{-})^{\ast})\bm{g}
\leq
\mathbb{1}.
\label{I-Tralpha1Abeta1pq-1-halpha1Abeta1pqg-1}
\end{equation}
are satisfied, and all solutions are given in parametric form by the relations
\begin{equation}
\bm{x}
=
(\alpha^{-1}\bm{A}
\oplus
\beta^{-1}\bm{p}\bm{q}^{-})^{\ast}
\bm{u},
\qquad
\bm{g}
\leq
\bm{u}
\leq
(\bm{h}^{-}
(\alpha^{-1}\bm{A}
\oplus
\beta^{-1}\bm{p}\bm{q}^{-})^{\ast})^{-}.
\label{E-xalpha1Abeta1pqu-guhalpha1Abeta1pq}
\end{equation}

\subsection{Determination of parameters}

To derive the Pareto frontier for the problem, we examine the conditions at \eqref{I-Tralpha1Abeta1pq-1-halpha1Abeta1pqg-1}. We consider the first condition and use identity \eqref{E-TrAB} for calculating
\begin{multline*}
\mathop\mathrm{Tr}(\alpha^{-1}\bm{A}\oplus\beta^{-1}\bm{p}\bm{q}^{-})
=
\bigoplus_{k=1}^{n}\alpha^{-k}\mathop\mathrm{tr}\bm{A}^{k}
\\
\oplus
\bigoplus_{k=1}^{n}
\beta^{-k}
\bigoplus_{m=0}^{n-k}
\bigoplus_{\substack{i_{1}+\cdots+i_{k}=m\\ i_{1},\ldots,i_{k}\geq0}}
\alpha^{-m}
(\bm{q}^{-}\bm{A}^{i_{1}}\bm{p}\cdots\bm{q}^{-}\bm{A}^{i_{k}}\bm{p}),
\end{multline*}
which allows to expand the condition as follows
\begin{equation*}
\bigoplus_{k=1}^{n}\alpha^{-k}\mathop\mathrm{tr}\bm{A}^{k}
\oplus
\bigoplus_{k=1}^{n}
\beta^{-k}
\bigoplus_{m=0}^{n-k}
\bigoplus_{\substack{i_{1}+\cdots+i_{k}=m\\ i_{1},\ldots,i_{k}\geq0}}
\alpha^{-m}
(\bm{q}^{-}\bm{A}^{i_{1}}\bm{p}\cdots\bm{q}^{-}\bm{A}^{i_{k}}\bm{p})
\leq
\mathbb{1}.
\end{equation*}

We replace the last inequality by the equivalent system of inequalities
\begin{align*}
\alpha^{-k}\mathop\mathrm{tr}\bm{A}^{k}
&\leq
\mathbb{1},
\\
\beta^{-k}
\bigoplus_{m=0}^{n-k}
\bigoplus_{\substack{i_{1}+\cdots+i_{k}=m\\ i_{1},\ldots,i_{k}\geq0}}
\alpha^{-m}
(\bm{q}^{-}\bm{A}^{i_{1}}\bm{p}\cdots\bm{q}^{-}\bm{A}^{i_{k}}\bm{p})
&\leq
\mathbb{1},
\qquad
k=1,\ldots,n.
\end{align*}

After rearranging terms to isolate powers of $\alpha$ and $\beta$ on the right-hand side and taking roots, we rewrite the system as
\begin{align*}
\mathop\mathrm{tr}\nolimits^{1/k}(\bm{A}^{k})
&\leq
\alpha,
\\
\bigoplus_{m=0}^{n-k}
\bigoplus_{\substack{i_{1}+\cdots+i_{k}=m\\ i_{1},\ldots,i_{k}\geq0}}
\alpha^{-m/k}
(\bm{q}^{-}\bm{A}^{i_{1}}\bm{p}\cdots\bm{q}^{-}\bm{A}^{i_{k}}\bm{p})^{1/k}
&\leq
\beta,
\qquad
k=1,\ldots,n.
\end{align*}

Aggregating the inequalities for $\alpha$ and then for $\beta$ yields the system
\begin{align*}
\alpha
&\geq
\bigoplus_{k=1}^{n}
\mathop\mathrm{tr}\nolimits^{1/k}(\bm{A}^{k})
=
\lambda,
\\
\beta
&\geq
\bigoplus_{k=1}^{n}
\bigoplus_{m=0}^{n-k}
\bigoplus_{\substack{i_{1}+\cdots+i_{k}=m\\ i_{1},\ldots,i_{k}\geq0}}
\alpha^{-m/k}
(\bm{q}^{-}\bm{A}^{i_{1}}\bm{p}\cdots\bm{q}^{-}\bm{A}^{i_{k}}\bm{p})^{1/k},
\end{align*}
where $\lambda$ denotes the spectral radius of the matrix $\bm{A}$.

To simplify the sum on the right-hand side of the second inequality, consider that part of this sum, which corresponds to $k=1$ and takes the form
\begin{equation*}
\bigoplus_{m=0}^{n-1}
\bigoplus_{i_{1}=m}
\alpha^{-m}
(\bm{q}^{-}\bm{A}^{i_{1}}\bm{p})
=
\bigoplus_{i=0}^{n-1}
\alpha^{-i}
(\bm{q}^{-}\bm{A}^{i}\bm{p}).
\end{equation*}

Let us verify that the rest of the sum is dominated by this part and thus can be eliminated. Indeed, with the condition that $i_{1}+\cdots+i_{k}=m$ where $1<k\leq n$ and $0\leq m\leq n-k$, we apply the tropical inequality between geometric and arithmetic means to obtain
\begin{multline*}
\alpha^{-m/k}
(\bm{q}^{-}\bm{A}^{i_{1}}\bm{p}\cdots\bm{q}^{-}\bm{A}^{i_{k}}\bm{p})^{1/k}
=
(\alpha^{-i_{1}}(\bm{q}^{-}\bm{A}^{i_{1}}\bm{p})\cdots\alpha^{-i_{k}}(\bm{q}^{-}\bm{A}^{i_{k}}\bm{p}))^{1/k}
\\\leq
\alpha^{-i_{1}}(\bm{q}^{-}\bm{A}^{i_{1}}\bm{p})\oplus\cdots\oplus\alpha^{-i_{k}}(\bm{q}^{-}\bm{A}^{i_{k}}\bm{p})
\leq
\bigoplus_{i=0}^{n-1}
\alpha^{-i}
(\bm{q}^{-}\bm{A}^{i}\bm{p}).
\end{multline*}

It follows from this inequality that the sum of all terms corresponding to $k>1$ satisfies the condition
\begin{equation*}
\bigoplus_{k=2}^{n}
\bigoplus_{m=0}^{n-k}
\bigoplus_{\substack{i_{1}+\cdots+i_{k}=m\\ i_{1},\ldots,i_{k}\geq0}}
\alpha^{-m/k}
(\bm{q}^{-}\bm{A}^{i_{1}}\bm{p}\cdots\bm{q}^{-}\bm{A}^{i_{k}}\bm{p})^{1/k}
\leq
\bigoplus_{i=0}^{n-1}
\alpha^{-i}
(\bm{q}^{-}\bm{A}^{i}\bm{p}),
\end{equation*}
and hence can be dropped without affecting the entire sum. 

As a result, we obtain the system of inequalities in the reduced form
\begin{equation}
\alpha
\geq
\bigoplus_{k=1}^{n}
\mathop\mathrm{tr}\nolimits^{1/k}(\bm{A}^{k}),
\qquad
\beta
\geq
\bigoplus_{k=0}^{n-1}
\alpha^{-k}
(\bm{q}^{-}\bm{A}^{k}\bm{p}).
\label{I-alphatr1kAk-I-betaAalphakqAkp}
\end{equation}

Next, we examine the second condition at \eqref{I-Tralpha1Abeta1pq-1-halpha1Abeta1pqg-1}. Application of identity \eqref{E-Apqast} to the Kleene star matrix, followed by the multiplication of the result by $\bm{h}^{-}$ on the left and by $\bm{g}$ on the right puts this condition into the form 
\begin{multline*}
\bm{h}^{-}\bm{g}
\oplus
\bigoplus_{k=1}^{n-1}
\alpha^{-k}(\bm{h}^{-}\bm{A}^{k}\bm{g})
\\\oplus
\bigoplus_{k=1}^{n-1}
\beta^{-k}
\bigoplus_{m=0}^{n-k-1}
\bigoplus_{\substack{i_{0}+i_{1}+\cdots+i_{k}=m\\i_{0},\ldots,i_{k}\geq0}}
\alpha^{-m}
(\bm{q}^{-}\bm{A}^{i_{1}}\bm{p}\cdots\bm{q}^{-}\bm{A}^{i_{k-1}}\bm{p})\bm{h}^{-}\bm{A}^{i_{0}}\bm{p}\bm{q}^{-}\bm{A}^{i_{k}}\bm{g}
\leq
\mathbb{1}.
\end{multline*}

We solve this inequality for $\alpha$ and $\beta$ in the same way as before. First note that the condition $\bm{g}\leq\bm{h}$ yields $\bm{h}^{-}\bm{g}\leq\mathbb{1}$, and hence the term $\bm{h}^{-}\bm{g}$ can be eliminated. Next, we take the other two terms to replace the inequality by a system of two inequalities.

After solving these inequalities and combining the solutions in the same way as above, we rewrite the system as
\begin{align*}
\alpha
&\geq
\bigoplus_{k=1}^{n-1}
(\bm{h}^{-}\bm{A}^{k}\bm{g})^{1/k},
\\
\beta
&\geq
\bigoplus_{m=0}^{n-2}
\bigoplus_{\substack{i+j=m\\i,j\geq0}}
\alpha^{-m}
(\bm{h}^{-}\bm{A}^{i}\bm{p})(\bm{q}^{-}\bm{A}^{j}\bm{g})
\\&\oplus
\bigoplus_{k=2}^{n-1}
\bigoplus_{m=0}^{n-k-1}
\bigoplus_{\substack{i_{0}+i_{1}+\cdots+i_{k}=m\\i_{0},\ldots,i_{k}\geq0}}
\alpha^{-m/k}
((\bm{q}^{-}\bm{A}^{i_{1}}\bm{p}\cdots\bm{q}^{-}\bm{A}^{i_{k-1}}\bm{p})\bm{h}^{-}\bm{A}^{i_{0}}\bm{p}\bm{q}^{-}\bm{A}^{i_{k}}\bm{g})^{1/k},
\end{align*}
where two sums on the right-hand side of the second inequality are obtained by separating all summands which correspond to $k=1$ from the others.

By coupling the last inequalities for $\alpha$ and $\beta$ with the corresponding inequalities from \eqref{I-alphatr1kAk-I-betaAalphakqAkp}, we form the system
\begin{align*}
\alpha
&\geq
\bigoplus_{k=1}^{n}
\mathop\mathrm{tr}\nolimits^{1/k}(\bm{A}^{k})
\oplus
\bigoplus_{k=1}^{n-1}
(\bm{h}^{-}\bm{A}^{k}\bm{g})^{1/k},
\\
\beta
&\geq
\bigoplus_{k=0}^{n-1}
\alpha^{-k}
(\bm{q}^{-}\bm{A}^{k}\bm{p})
\oplus
\bigoplus_{m=0}^{n-2}
\bigoplus_{\substack{i+j=m\\i,j\geq0}}
\alpha^{-m}
(\bm{h}^{-}\bm{A}^{i}\bm{p})(\bm{q}^{-}\bm{A}^{j}\bm{g})
\\&\oplus
\bigoplus_{k=2}^{n-1}
\bigoplus_{m=0}^{n-k-1}
\bigoplus_{\substack{i_{0}+i_{1}+\cdots+i_{k}=m\\i_{0},\ldots,i_{k}\geq0}}
\alpha^{-m/k}
((\bm{q}^{-}\bm{A}^{i_{1}}\bm{p}\cdots\bm{q}^{-}\bm{A}^{i_{k-1}}\bm{p})\bm{h}^{-}\bm{A}^{i_{0}}\bm{p}\bm{q}^{-}\bm{A}^{i_{k}}\bm{g})^{1/k}.
\end{align*}

To simplify the right-hand side of the second inequality, we verify that the first two sums on this side dominate the third sum. By using the inequality of geometric and arithmetic means once again, we obtain for all $k=2,\ldots,n-1$, $m=0,\ldots,n-k-1$ and $i_{0}+i_{1}+\cdots+i_{k}=m$, the inequalities
\begin{multline*}
\alpha^{-m/k}
((\bm{q}^{-}\bm{A}^{i_{1}}\bm{p}\cdots\bm{q}^{-}\bm{A}^{i_{k-1}}\bm{p})\bm{h}^{-}\bm{A}^{i_{0}}\bm{p}\bm{q}^{-}\bm{A}^{i_{k}}\bm{g})^{1/k}
\\\leq
(\alpha^{-i_{1}}(\bm{q}^{-}\bm{A}^{i_{1}}\bm{p})\oplus\cdots\oplus\alpha^{-i_{k-1}}(\bm{q}^{-}\bm{A}^{i_{k-1}}\bm{p}))\oplus\alpha^{-(i_{0}+i_{k})}(\bm{h}^{-}\bm{A}^{i_{0}}\bm{p}\bm{q}^{-}\bm{A}^{i_{k}}\bm{g})
\\\leq
\bigoplus_{i=0}^{n-1}\alpha^{-i}(\bm{q}^{-}\bm{A}^{i}\bm{p})
\oplus
\bigoplus_{l=0}^{n-2}
\bigoplus_{\substack{i+j=l\\i,j\geq0}}
\alpha^{-l}
(\bm{h}^{-}\bm{A}^{i}\bm{p})(\bm{q}^{-}\bm{A}^{j}\bm{g}),
\end{multline*}
which show that each summand of the third sum is not greater than the first two sums, and hence the third sum can be eliminated.

As a result, the system of inequalities for $\alpha$ and $\beta$ takes the form
\begin{equation}
\begin{aligned}
\alpha
&\geq
\bigoplus_{k=1}^{n}
\mathop\mathrm{tr}\nolimits^{1/k}(\bm{A}^{k})
\oplus
\bigoplus_{k=1}^{n-1}
(\bm{h}^{-}\bm{A}^{k}\bm{g})^{1/k},
\\
\beta
&\geq
\bm{q}^{-}\bm{p}
\oplus
\bm{h}^{-}\bm{p}\bm{q}^{-}\bm{g}
\oplus
\bigoplus_{k=1}^{n-1}\alpha^{-k}
(\bm{q}^{-}\bm{A}^{k}\bm{p})
\oplus
\bigoplus_{k=1}^{n-2}
\alpha^{-k}
\bigoplus_{\substack{i+j=k\\i,j\geq0}}
(\bm{h}^{-}\bm{A}^{i}\bm{p})(\bm{q}^{-}\bm{A}^{j}\bm{g}).
\end{aligned}
\label{I-alpha-tr1kAkhAkg1k-beta-qp}
\end{equation}

\subsection{Derivation of Pareto frontier}

We now examine the system at \eqref{I-alpha-tr1kAkhAkg1k-beta-qp} to derive a representation for the Pareto frontier in terms of the parameters $\alpha$ and $\beta$. To simplify further formulas, we use the notation
\begin{equation}
\lambda
=
\bigoplus_{k=1}^{n}
\mathop\mathrm{tr}\nolimits^{1/k}(\bm{A}^{k}),
\qquad
\mu
=
\bigoplus_{k=1}^{n-1}
(\bm{h}^{-}\bm{A}^{k}\bm{g})^{1/k},
\qquad
\nu
=
\bm{q}^{-}\bm{p}
\oplus
\bm{h}^{-}\bm{p}\bm{q}^{-}\bm{g},
\label{E-lambda-mu-nu}
\end{equation}
and introduce the functions
\begin{equation}
\begin{aligned}
G(s)
&=
\bigoplus_{k=1}^{n-1}
s^{-k}
(\bm{q}^{-}\bm{A}^{k}\bm{p})
\oplus
\bigoplus_{k=1}^{n-2}
s^{-k}
\bigoplus_{\substack{i+j=k\\i,j\geq0}}
(\bm{h}^{-}\bm{A}^{i}\bm{p})(\bm{q}^{-}\bm{A}^{j}\bm{g}),
\quad
s>\mathbb{0};
\\
H(t)
&=
\bigoplus_{k=1}^{n-1}
t^{-1/k}
(\bm{q}^{-}\bm{A}^{k}\bm{p})^{1/k}
\oplus
\bigoplus_{k=1}^{n-2}
t^{-1/k}
\bigoplus_{\substack{i+j=k\\i,j\geq0}}
(\bm{h}^{-}\bm{A}^{i}\bm{p})^{1/k}(\bm{q}^{-}\bm{A}^{j}\bm{g})^{1/k},
\quad
t>\mathbb{0}.
\end{aligned}
\label{E-G-H}
\end{equation}

We note that both functions are monotone decreasing. Furthermore, it is not difficult to verify by direct calculation that the inequalities
\begin{equation*}
G(s)
\leq
t,
\qquad
H(t)
\leq
s
\end{equation*}
are equivalent in the sense that all solutions of the first inequality with respect to $s$ are represented by the second inequality and vice versa (see, also, \cite{Krivulin2020Using}).

Finally, with the new notation, we represent the system \eqref{I-alpha-tr1kAkhAkg1k-beta-qp} in the form 
\begin{equation}
\alpha
\geq
\lambda\oplus\mu,
\qquad
\beta
\geq
\nu
\oplus
G(\alpha).
\label{I-alpha-lambdamu-beta-muGalpha}
\end{equation}

Consider the area of points $(\alpha,\beta)$, which is given by the system of inequalities at \eqref{I-alpha-lambdamu-beta-muGalpha}. To determine the Pareto frontier as the image of the set of Pareto optimal solutions, we need to examine the boundary of this area because each interior point apparently corresponds to a solution that is dominated and hence is not Pareto-optimal.

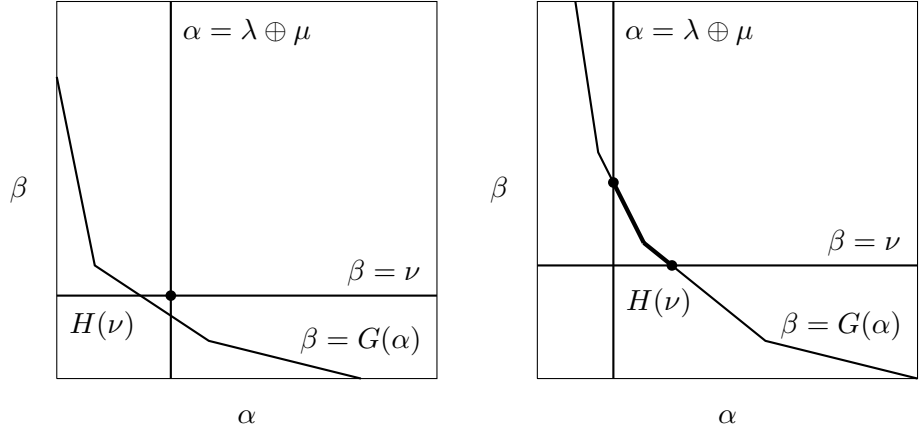
\begin{figure}[ht]
\begin{tikzpicture}

\draw (0.5,0.5) -- (5.5,0.5) -- (5.5,5.5) -- (0.5,5.5) -- (0.5,0.5);
\node at (0.0,3.0) {$\beta$};
\node at (3.0,0.0) {$\alpha$};

\draw [thick] (2.0,0.5) -- (2.0,5.5);
\draw [thick] (0.5,1.6) -- (5.5,1.6);

\draw [thick] (0.5,4.5) -- (1.0, 2.0) -- (2.5,1.0) -- (4.5,0.5);

\fill (2.0,1.6) circle (2pt);

\node at (3.0,5.1) {$\alpha=\lambda\oplus\mu$};

\node at (4.8,1.9) {$\beta=\nu$};

\node at (1.1,1.2) {$H(\nu)$};

\node at (4.5,1.0) {$\beta=G(\alpha)$};

\end{tikzpicture}
\hspace{3mm}
\begin{tikzpicture}

\draw (0.5,0.5) -- (5.5,0.5) -- (5.5,5.5) -- (0.5,5.5) -- (0.5,0.5);
\node at (0.0,3.0) {$\beta$};
\node at (3.0,0.0) {$\alpha$};

\draw [thick] (1.5,0.5) -- (1.5,5.5);
\draw [thick] (0.5,2.0) -- (5.5,2.0);

%\draw [thick] (1.0,5.5) -- (1.3,3.5) -- (2.0,2.3) -- (3.5,1.0) -- (5.5,0.5);
\draw [thick] (1.0,5.5) -- (1.3,3.5) -- (1.9,2.3) -- (3.5,1.0) -- (5.5,0.5);

\fill (1.5,3.1) circle (2pt);
\fill (2.27,2.0) circle (2pt);

%\fill (1.9,2.3) circle (2pt);

\draw [ultra thick] (1.5,3.1) -- (1.9,2.3);
\draw [ultra thick] (1.9,2.3) -- (2.27,2.0);

\node at (2.5,5.1) {$\alpha=\lambda\oplus\mu$};

\node at (4.8,2.3) {$\beta=\nu$};

\node at (2.1,1.5) {$H(\nu)$};

\node at (4.5,1.2) {$\beta=G(\alpha)$};
\end{tikzpicture}
\caption{Examples of Pareto frontiers in terms of $\mathbb{R}_{\max,+}$ in the form of a point (left) and a segment (right).}
\label{F-EPF}
\end{figure}

The area is bounded from the left by the vertical line $\alpha=\lambda\oplus\mu$, from the lower left by the curve $\beta=\nu\oplus G(\alpha)$ and from below by the horizontal line $\beta=\nu$. Since the function $G(\alpha)$ is monotone decreasing, the Pareto frontier for the problem is a segment of the curve that lies right of the vertical and above the horizontal lines, or a single point if the curve lies below the intersection of these lines, as shown in Fig.~\ref{F-EPF}. 

To describe the Pareto frontier and related Pareto-optimal solutions, we examine two cases. First, we assume that the following condition holds:
\begin{equation*}
\lambda\oplus\mu
\geq
H(\nu).
\end{equation*}

Under this condition, it follows from the inequality $\alpha\geq\lambda\oplus\nu$ that the inequality $\alpha\geq H(\nu)$ is satisfied. Solving the last inequality with respect to $\nu$ yields $\nu\geq G(\alpha)$. As a result, system \eqref{I-alpha-lambdamu-beta-muGalpha} becomes
\begin{equation*}
\alpha
\geq
\lambda\oplus\mu,
\qquad
\beta
\geq
\nu,
\end{equation*}
which gives the Pareto frontier reduced to a single point $(\alpha,\beta)$ where
\begin{equation*}
\alpha
=
\lambda\oplus\mu
=
\bigoplus_{k=1}^{n}
\mathop\mathrm{tr}\nolimits^{1/k}(\bm{A}^{k})
\oplus
\bigoplus_{k=1}^{n-1}
(\bm{h}^{-}\bm{A}^{k}\bm{g})^{1/k},
\qquad
\beta
=
\nu
=
\bm{q}^{-}\bm{p}
\oplus
\bm{h}^{-}\bm{p}\bm{q}^{-}\bm{g}.
\end{equation*}

An example of the Pareto frontier in terms of the semifield $\mathbb{R}_{\max,+}$ (max-plus algebra) for this case is shown in Fig.~\ref{F-EPF} (left) by the thick dot in the intersection of the lines $\alpha=\lambda\oplus\mu$ and $\beta=\nu$. 

Now suppose that the opposite condition is valid in the form
\begin{equation*}
\lambda\oplus\mu
<
H(\nu).
\end{equation*}

Then, system \eqref{I-alpha-lambdamu-beta-muGalpha} defines an area which is given by the conditions
\begin{equation*}
\lambda\oplus\mu
\leq
\alpha
\leq
H(\nu),
\qquad
\beta
\geq
\nu
\oplus
G(\alpha),
\end{equation*}
with its lower left boundary defined by $\lambda\oplus\mu\leq\alpha\leq H(\nu)$ and $\beta=\nu\oplus G(\alpha)$.

Let us verify that, under the condition $\lambda\oplus\mu\leq\alpha\leq H(\nu)$, we have $G(\alpha)\geq\nu$, and hence the last equality can be reduced to $\beta=G(\alpha)$.

Since $G(\alpha)$ is a monotone decreasing function of $\alpha$, it is sufficient to verify that $G(\alpha)\geq\nu$ for $\alpha=H(\nu)$. We consider the equality
\begin{multline*}
H(\nu)
=
\bigoplus_{k=1}^{n-1}
(\bm{q}^{-}\bm{p}
\oplus
\bm{h}^{-}\bm{p}\bm{q}^{-}\bm{g})^{-1/k}
(\bm{q}^{-}\bm{A}^{k}\bm{p})^{1/k}
\\\oplus
\bigoplus_{k=1}^{n-2}
(\bm{q}^{-}\bm{p}
\oplus
\bm{h}^{-}\bm{p}\bm{q}^{-}\bm{g})^{-1/k}
\bigoplus_{\substack{i+j=k\\i,j\geq0}}
(\bm{h}^{-}\bm{A}^{i}\bm{p})^{1/k}(\bm{q}^{-}\bm{A}^{j}\bm{g})^{1/k},
\end{multline*}
and note that this equality is valid in two cases.

For the first case, we assume that there exists an index $m$ such that
%\begin{multline*}
\begin{equation*}
H(\nu)
=
%\bigoplus_{k=1}^{n-1}
%(\bm{q}^{-}\bm{p}
%\oplus
%\bm{h}^{-}\bm{p}\bm{q}^{-}\bm{g})^{-1/k}
%(\bm{q}^{-}\bm{A}^{k}\bm{p})^{1/k}
%\\=
(\bm{q}^{-}\bm{p}
\oplus
\bm{h}^{-}\bm{p}\bm{q}^{-}\bm{g})^{-1/m}
(\bm{q}^{-}\bm{A}^{m}\bm{p})^{1/m}.
\end{equation*}
%\end{multline*}
 
Then, under the assumption that $\alpha=H(\nu)$, we obtain
\begin{multline*}
G(\alpha)
=
\bigoplus_{k=1}^{n-1}
\alpha^{-k}
(\bm{q}^{-}\bm{A}^{k}\bm{p})
\oplus
\bigoplus_{k=1}^{n-2}
\alpha^{-k}
\bigoplus_{\substack{i+j=k\\i,j\geq0}}
(\bm{h}^{-}\bm{A}^{i}\bm{p})(\bm{q}^{-}\bm{A}^{j}\bm{g})
\\\geq
\bigoplus_{k=1}^{n-1}
\alpha^{-k}
(\bm{q}^{-}\bm{A}^{k}\bm{p})
\geq
\alpha^{-m}
(\bm{q}^{-}\bm{A}^{m}\bm{p})
=
\bm{q}^{-}\bm{p}
\oplus
\bm{h}^{-}\bm{p}\bm{q}^{-}\bm{g}
=
\nu.
\end{multline*}

The other case involves the condition that for some $m$ the following equality holds:
%\begin{multline*}
\begin{equation*}
H(\nu)
=
%\bigoplus_{k=1}^{n-2}
%(\bm{q}^{-}\bm{p}
%\oplus
%\bm{h}^{-}\bm{p}\bm{q}^{-}\bm{g})^{-1/k}
%\bigoplus_{\substack{i+j=k\\i,j\geq0}}
%(\bm{h}^{-}\bm{A}^{i}\bm{p})^{1/k}(\bm{q}^{-}\bm{A}^{j}\bm{g})^{1/k}
%\\=
(\bm{q}^{-}\bm{p}
\oplus
\bm{h}^{-}\bm{p}\bm{q}^{-}\bm{g})^{-1/m}
\bigoplus_{\substack{i+j=k\\i,j\geq0}}
(\bm{h}^{-}\bm{A}^{i}\bm{p})^{1/m}(\bm{q}^{-}\bm{A}^{j}\bm{g})^{1/m},
%\end{multline*}
\end{equation*}
which with setting $\alpha=H(\nu)$ yields
\begin{multline*}
G(\alpha)
\geq
\bigoplus_{k=1}^{n-2}
\alpha^{-k}
\bigoplus_{\substack{i+j=k\\i,j\geq0}}
(\bm{h}^{-}\bm{A}^{i}\bm{p})(\bm{q}^{-}\bm{A}^{j}\bm{g})
\\\geq
\alpha^{-m}
\bigoplus_{\substack{i+j=m\\i,j\geq0}}
(\bm{h}^{-}\bm{A}^{i}\bm{p})(\bm{q}^{-}\bm{A}^{j}\bm{g})
=
\bm{q}^{-}\bm{p}
\oplus
\bm{h}^{-}\bm{p}\bm{q}^{-}\bm{g}
=
\nu.
\end{multline*}

Since $G(\alpha)\geq\nu$ in both cases, the description of the Pareto frontier can be reduced to the system
\begin{equation*}
\lambda\oplus\mu
\leq
\alpha
\leq
H(\nu),
\qquad
\beta
=
G(\alpha).
\end{equation*}

This case is illustrated in Fig.~\ref{F-EPF} (right) where the Pareto frontier is depicted by a thick segment that is cut by the vertical line $\alpha=\lambda\oplus\mu$ and horizontal line $\beta=\nu$ from the curve $\beta=G(\alpha)$.

\subsection{Pareto-optimal solution}

We are now in a position to summarize the obtained results by the next statement.
\begin{theorem}
\label{T-minxAxxpqx-gxh}
For any nonzero matrix $\bm{A}$, nonzero vectors $\bm{p}$, $\bm{q}$, $\bm{g}$ and $\bm{h}$ such that $\bm{h}^{-}\bm{g}\leq\mathbb{1}$, with the notation \eqref{E-lambda-mu-nu}--\eqref{E-G-H} the following statements hold:
\begin{enumerate}
\item
If $\lambda\oplus\mu\geq H(\nu)$, the Pareto frontier for problem \eqref{P-minxAxxpqx-gxh} degenerates into the single point $(\alpha,\beta)$ with $\alpha=\lambda\oplus\mu$ and $\beta=\nu$.
\item
Otherwise, the Pareto frontier is the segment given by the conditions
\begin{equation*}
\lambda\oplus\mu
\leq
\alpha
\leq
H(\nu),
\qquad
\beta
=
G(\alpha).
\end{equation*}
\end{enumerate}

All Pareto-optimal solutions are represented in parametric form as
\begin{equation*}
\bm{x}
=
(\alpha^{-1}\bm{A}
\oplus
\beta^{-1}\bm{p}\bm{q}^{-})^{\ast}
\bm{u},
\qquad
\bm{g}
\leq
\bm{u}
\leq
(\bm{h}^{-}
(\alpha^{-1}\bm{A}
\oplus
\beta^{-1}\bm{p}\bm{q}^{-})^{\ast})^{-}.
\end{equation*}
\end{theorem}

We now briefly discuss the computational complexity of the solution obtained. First, we note that the computational time required to obtain a solution vector depends on the time spent on calculating values of $\lambda$, $\mu$, $G$ and $H$, and evaluating the matrix $(\alpha^{-1}\bm{A}\oplus\beta^{-1}\bm{p}\bm{q}^{-})^{\ast}$. As it easy to see, the most computationally demanding task involved is obtaining the first $n$ powers of the matrices $\bm{A}$ and $\alpha^{-1}\bm{A}\oplus\beta^{-1}\bm{p}\bm{q}^{-}$. Since direct multiplication of two matrices of order $n$ takes at most $O(n^{3})$ operations, the time to obtain these powers can be estimated $O(n^{4})$. Similar reasoning shows that both functions $G$ and $H$ can be computed in the same time. As a result, the overall computational complexity of the solution is no more than $O(n^{4})$.

We also note that the compact vector form of the solution obtained, which involves a fixed number of simple matrix and vector operations, makes the result easily scalable to high-dimensional problems and offers a strong potential for efficient implementation on parallel computing platforms.

\section{Application to bi-objective project scheduling problems}
\label{S-ABOSP}

Consider the bi-objective project scheduling problem at \eqref{P-minyixi-yixi-aijxjyi-gixihi}, and note that the representation of both the objective functions and the constraints involves only the arithmetic operations of maximum, addition and additive inversion (subtraction). As a consequence, we can rewrite \eqref{P-minyixi-yixi-aijxjyi-gixihi} in the max-plus algebra setting by changing the operation symbols to obtain a tropical optimization problem that is to
\begin{equation*}
\begin{aligned}
&
\text{minimize}
&&
\left\{
\bigoplus_{i=1}^{n}y_{i}x_{i}^{-1},\ 
\bigoplus_{i=1}^{n}y_{i}\bigoplus_{j=1}^{n}x_{j}^{-1}
\right\};
\\
&
\text{subject to}
&&
\bigoplus_{j=1}^{n}a_{ij}x_{j}
=
y_{i},
\\
&&&
g_{i}
\leq
x_{i}
\leq
h_{i},
\qquad
i=1,\ldots,n.
%\label{P-minyixi-yixi-aijxjyi-gixihi}
\end{aligned}
\end{equation*}

Furthermore, we introduce the following matrix and vectors:
\begin{equation*}
\bm{A}
=
(a_{ij}),
\qquad
\bm{x}
=
(x_{i}),
\qquad
\bm{y}
=
(y_{i}),
\qquad
\bm{g}
=
(g_{i}),
\qquad
\bm{h}
=
(h_{i}),
\end{equation*}
and represent the problem in the vector form 
\begin{equation*}
\begin{aligned}
&
\text{minimize}
&&
\left\{
\bm{x}^{-}\bm{y},\ 
\bm{x}^{-}\bm{1}\bm{1}^{T}\bm{y}
\right\};
\\
&
\text{subject to}
&&
\bm{A}\bm{x}
=
\bm{y},
\quad
\bm{g}
\leq
\bm{x}
\leq
\bm{h}.
\end{aligned}
\end{equation*}

After substitution $\bm{y}=\bm{A}\bm{x}$ into the objective functions, we obtain a problem in the form of \eqref{P-minxAxxpqx-gxh} where $\bm{p}=\bm{1}$ and $\bm{q}^{-}=\bm{1}^{T}\bm{A}$. We now exploit Theorem~\ref{T-minxAxxpqx-gxh} to describe a complete Pareto-optimal solution to the problem. 

To rewrite the statement of Theorem~\ref{T-minxAxxpqx-gxh} in terms of the problem under study, we first observe that $\bm{x}^{-}\bm{1}=\|\bm{x}^{-}\|$ and $\bm{1}^{T}\bm{A}\bm{x}=\|\bm{A}\bm{x}\|$. Then, problem \eqref{P-minxAxxpqx-gxh} takes the form
\begin{equation}
\begin{aligned}
&
\text{minimize}
&&
\{\bm{x}^{-}\bm{A}\bm{x},\ 
\|\bm{x}^{-}\|\|\bm{A}\bm{x}\|\};
\\
&
\text{subject to}
&&
\bm{g}
\leq
\bm{x}
\leq
\bm{h}.
\end{aligned}
\label{P-minxAxxAx-gleqxleqh}
\end{equation}

Furthermore, we refine the notation \eqref{E-lambda-mu-nu} as
\begin{equation*}
\lambda
=
\bigoplus_{k=1}^{n}
\mathop\mathrm{tr}\nolimits^{1/k}(\bm{A}^{k}),
\qquad
\mu
=
\bigoplus_{k=1}^{n-1}
(\bm{h}^{-}\bm{A}^{k}\bm{g})^{1/k},
\qquad
\nu
=
\|\bm{A}\|\oplus\|\bm{h}^{-}\|\|\bm{A}\bm{g}\|,
\end{equation*}
and the notation \eqref{E-G-H} as
\begin{equation*}
\begin{aligned}
G(s)
&=
\bigoplus_{k=1}^{n-1}
s^{-k}
\|\bm{A}^{k+1}\|
\oplus
\bigoplus_{k=1}^{n-2}
s^{-k}
\bigoplus_{\substack{i+j=k\\i,j\geq0}}
\|\bm{h}^{-}\bm{A}^{i}\|\|\bm{A}^{j+1}\bm{g}\|,
\\
H(t)
&=
\bigoplus_{k=1}^{n-1}
t^{-1/k}
\|\bm{A}^{k+1}\|^{1/k}
\oplus
\bigoplus_{k=1}^{n-2}
t^{-1/k}
\bigoplus_{\substack{i+j=k\\i,j\geq0}}
\|\bm{h}^{-}\bm{A}^{i}\|^{1/k}\|\bm{A}^{j+1}\bm{g}\|^{1/k}.
\end{aligned}
\end{equation*}

The solution given by Theorem~\ref{T-minxAxxpqx-gxh} turns into the next result.
\begin{corollary}
\label{C-minxAxxAx-gleqxleqh}
For any matrix $\bm{A}$ and vectors $\bm{g}$ and $\bm{h}$ such that $\bm{h}^{-}\bm{g}\leq\mathbb{1}$, the following statements hold:
\begin{enumerate}
\item
If $\lambda\oplus\mu\geq H(\nu)$, the Pareto frontier for problem \eqref{P-minxAxxAx-gleqxleqh} degenerates into the single point $(\alpha,\beta)$ with $\alpha=\lambda\oplus\mu$ and $\beta=\nu$.
\item
Otherwise, the Pareto frontier is the segment given by the conditions
\begin{equation*}
\lambda\oplus\mu
\leq
\alpha
\leq
H(\nu),
\qquad
\beta
=
G(\alpha).
\end{equation*}
\end{enumerate}

All Pareto-optimal solutions are represented in parametric form as
\begin{equation*}
\bm{x}
=
(\alpha^{-1}\bm{A}
\oplus
\beta^{-1}\bm{1}\bm{1}^{T})^{\ast}
\bm{u},
\qquad
\bm{g}
\leq
\bm{u}
\leq
(\bm{h}^{-}
(\alpha^{-1}\bm{A}
\oplus
\beta^{-1}\bm{1}\bm{1}^{T})^{\ast})^{-}.
\end{equation*}
\end{corollary}
 
\begin{figure}[ht]
\hspace{-2mm}
\begin{tikzpicture}

\draw (0.5,0.5) -- (5.5,0.5) -- (5.5,5.5) -- (0.5,5.5) -- (0.5,0.5);

%\node at (0.0,5.5) {$\beta$};
\node at (0.1,5.5) {$\beta$};
\draw (0.5,1.5) -- (0.6,1.5);
\draw (0.5,2.5) -- (0.6,2.5);
\draw (0.5,3.5) -- (0.6,3.5);
\draw (0.5,4.5) -- (0.6,4.5);
\node at (0.2,0.2) {$0$};
\node at (0.2,1.6) {$1$};
\node at (0.2,2.6) {$2$};
\node at (0.2,3.6) {$3$};
\node at (0.2,4.6) {$4$};

\node at (5.5,0.0) {$\alpha$};
\draw (1.5,0.5) -- (1.5,0.6);
\draw (2.5,0.5) -- (2.5,0.6);
\draw (3.5,0.5) -- (3.5,0.6);
\draw (4.5,0.5) -- (4.5,0.6);
\node at (1.5,0.2) {$1$};
\node at (2.5,0.2) {$2$};
\node at (3.5,0.2) {$3$};
\node at (4.5,0.2) {$4$};

\draw [thick] (2.0,0.5) -- (2.0,5.5);
\draw [thick] (0.5,2.5) -- (5.5,2.5);

\draw [thick] (0.5,5.5) -- (1.5,3.5);
\draw [thick] (1.5,3.5) -- (4.5,0.5);

\fill (2.0,3) circle (2pt);
\fill (2.5,2.5) circle (2pt);

\draw [ultra thick] (2.0,3) -- (2.5,2.5);

\node at (2.9,5.1) {$\alpha=3/2$};

\node at (4.8,2.8) {$\beta=2$};

\node at (2.5,1.7) {$H(2)$};

\node at (4.5,1.7) {$\beta=G(\alpha)$};
\end{tikzpicture}
\hspace{0mm}
\newcommand\transparentcube[3]{%

\filldraw [gray!80!white] (tpp cs:x=0,y=0,z=0) -- (tpp cs:x=#1,y=#2,z=#3) -- (tpp cs:x=#1,y=0,z=0) -- cycle; 
\draw [very thick] (tpp cs:x=0,y=0,z=0) -- (tpp cs:x=#1,y=#2,z=#3) -- (tpp cs:x=#1,y=0,z=0) -- cycle; 

%\fill (tpp cs:x=0,y=0,z=0) circle (2pt);
%\fill (tpp cs:x=#1,y=#2,z=#3) circle (2pt);
%\fill (tpp cs:x=#1,y=0,z=0) circle (2pt);

\filldraw [gray] (tpp cs:x=2*#1/3,y=#2/3,z=0) -- (tpp cs:x=#1,y=2*#2/3,z=#3/3) -- (tpp cs:x=#1,y=#2/3,z=0) -- cycle; 
\draw [very thick] (tpp cs:x=2*#1/3,y=#2/3,z=0) -- (tpp cs:x=#1,y=2*#2/3,z=#3/3) -- (tpp cs:x=#1,y=#2/3,z=0) -- cycle; 

%\fill (tpp cs:x=2*#1/3,y=#2/3,z=0) circle (2pt);
%\fill (tpp cs:x=#1,y=2*#2/3,z=#3/3) circle (2pt);
%\fill (tpp cs:x=#1,y=#2/3,z=0) circle (2pt);

\draw [dashed] (tpp cs:x=0,y=0,z=0) -- (tpp cs:x=#1,y=#2/2,z=0);
\draw [dashed] (tpp cs:x=#1,y=#2/2,z=0) -- (tpp cs:x=#1,y=#2,z=#3);

\fill (tpp cs:x=#1,y=#2/2,z=0) circle (2pt);

\draw (tpp cs:x=0,y=0,z=0) -- (tpp cs:x=0,y=0,z=#3);
\draw (tpp cs:x=0,y=0,z=0) -- (tpp cs:x=0,y=#2,z=0);
\draw (tpp cs:x=0,y=0,z=0) -- (tpp cs:x=#1,y=0,z=0);

\draw (tpp cs:x=0,y=0,z=#3) -- (tpp cs:x=0,y=#2,z=#3);
\draw (tpp cs:x=0,y=#2,z=0) -- (tpp cs:x=#1,y=#2,z=0);
\draw (tpp cs:x=#1,y=0,z=0) -- (tpp cs:x=#1,y=0,z=#3);

\draw (tpp cs:x=#1,y=#2,z=#3) -- (tpp cs:x=#1,y=#2,z=0);
\draw (tpp cs:x=#1,y=#2,z=#3) -- (tpp cs:x=#1,y=0,z=#3);
\draw (tpp cs:x=#1,y=#2,z=#3) -- (tpp cs:x=0,y=#2,z=#3);

\draw (tpp cs:x=#1,y=#2,z=0) -- (tpp cs:x=#1,y=0,z=0);
\draw (tpp cs:x=#1,y=0,z=#3) -- (tpp cs:x=0,y=0,z=#3);
\draw (tpp cs:x=0,y=#2,z=#3) -- (tpp cs:x=0,y=#2,z=0);

\node at (tpp cs:x=#1/2,y=-#2/8,z=-#1/18) {$x_{1}$};
\node at (tpp cs:x=#1,y=-#2/7,z=0) {$1$};

\node at (tpp cs:x=#1+#1/8,y=#2/2+#2/12,z=-#1/8) {$x_{2}$};
\node at (tpp cs:x=#1,y=#2/8,z=-#3/12) {$0$};
\node at (tpp cs:x=#1+#1/8,y=#2,z=-#3/18) {$1$};

\node at (tpp cs:x=0,y=-#1/8,z=#3/2) {$x_{3}$};
\node at (tpp cs:x=0,y=-#2/12,z=0) {$0$};
\node at (tpp cs:x=0,y=-#2/12,z=#3+#3/18) {$1$};
}
\begin{tikzpicture}[3d view={60}{15}]
\transparentcube{4}{4}{4}
\end{tikzpicture}
\caption{The Pareto frontier (left) and solutions sets (right) in Example~\ref{X-P3}.}
\label{F-PFSS}
\end{figure}
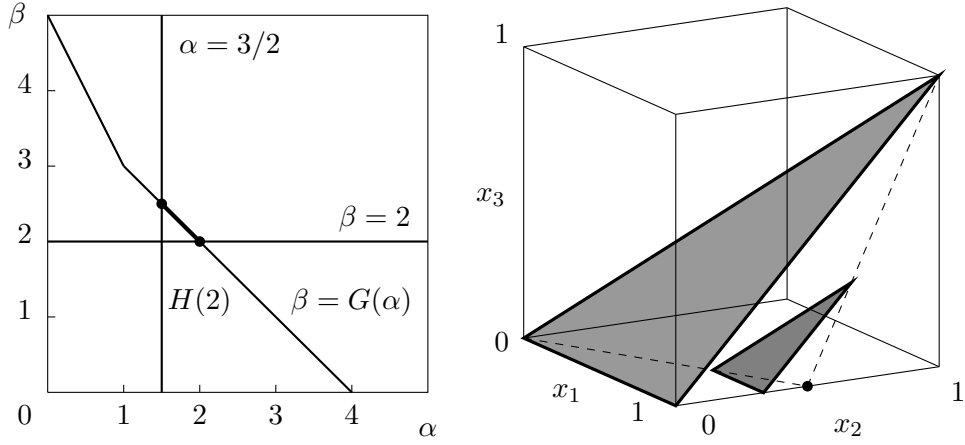

To demonstrate application of the result obtained, we offer an example of solution of a three-dimensional problem, which, in particular, shows that the Pareto frontier may be a segment rather than a single point. Although the example deals with a somewhat artificial problem, it clearly demonstrates the proposed computational technique which involves a fixed number of simple matrix-vector operations and hence allows natural scalability to handle problems of high dimension. 
\begin{example}
\label{X-P3}
Consider a project that performs $n=3$ activities under start-finish, release start and release end time constraints given by the following matrix and vectors:
\begin{equation*}
\bm{A}
=
\left(
\begin{array}{ccc}
1 & 2 & 2
\\
1 & 1 & 2
\\
\mathbb{0} & 0 & 1
\end{array}
\right),
\qquad
\bm{g}
=
\left(
\begin{array}{c}
0
\\
0
\\
0
\end{array}
\right),
\qquad
\bm{h}
=
\left(
\begin{array}{c}
1
\\
2
\\
2
\end{array}
\right).
\end{equation*}

To find an optimal schedule by minimizing both the maximum flow-time of activities and the project makespan, we apply Corollary~\ref{C-minxAxxAx-gleqxleqh} in the max-plus algebra setting. By convention, we represent the numerical constants (including negative integer and rational numbers) in the ordinary notation, whereas all algebraic operations are considered in terms of max-plus algebra.

First, we form the matrices
\begin{equation*}
\bm{A}^{2}
=
\left(
\begin{array}{ccc}
3 & 3 & 4
\\
2 & 3 & 3
\\
1 & 1 & 2
\end{array}
\right),
\qquad
\bm{A}^{3}
=
\left(
\begin{array}{ccc}
4 & 5 & 5
\\
4 & 4 & 5
\\
2 & 3 & 3
\end{array}
\right),
\qquad
\bm{1}\bm{1}^{T}\bm{A}
=
\left(
\begin{array}{ccc}
1 & 2 & 2
\\
1 & 2 & 2
\\
1 & 2 & 2
\end{array}
\right),
\end{equation*}
and then obtain
\begin{equation*}
\mathop\mathrm{tr}\bm{A}
=
1,
\quad
\mathop\mathrm{tr}\bm{A}^{2}
=
3,
\quad
\mathop\mathrm{tr}\bm{A}^{3}
=
4,
\quad
\|\bm{A}\|
=
2,
\quad
\|\bm{A}^{2}\|
=
4,
\quad
\|\bm{A}^{3}\|
=
5.
\end{equation*}

Next, we calculate the vectors
\begin{equation*}
\bm{h}^{-}
=
\left(
\begin{array}{rrr}
-1 & -2 & -2
\end{array}
\right),
\qquad
\bm{h}^{-}\bm{A}
=
\left(
\begin{array}{ccc}
0 & 1 & 1
\end{array}
\right),
\qquad
\bm{A}\bm{g}
=
\left(
\begin{array}{c}
2
\\
2
\\
1
\end{array}
\right),
\qquad
\bm{A}^{2}\bm{g}
=
\left(
\begin{array}{c}
4
\\
3
\\
2
\end{array}
\right),
\end{equation*}
which allow us to find the scalars
\begin{gather*}
\|\bm{h}^{-}\|
=
-1,
\qquad
\|\bm{h}^{-}\bm{A}\|
=
1,
\qquad
\bm{h}^{-}\bm{A}\bm{g}
=
1,
\qquad
\bm{h}^{-}\bm{A}^{2}\bm{g}
=
3,
\\
\|\bm{A}\bm{g}\|
=
2,
\qquad
\|\bm{A}^{2}\bm{g}\|
=
4.
\end{gather*}

Using the above results yields
\begin{gather*}
\lambda
=
\mathop\mathrm{tr}\bm{A}
\oplus
\mathop\mathrm{tr}\nolimits^{1/2}(\bm{A}^{2})
\oplus
\mathop\mathrm{tr}\nolimits^{1/3}(\bm{A}^{3})
=
3/2,
\qquad
\mu
=
\bm{h}^{-}\bm{A}\bm{g}
\oplus
(\bm{h}^{-}\bm{A}^{2}\bm{g})^{1/2}
=
3/2,
\\
\nu
=
\|\bm{A}\|\oplus\|\bm{h}^{-}\|\|\bm{A}\bm{g}\|
=
2.
\end{gather*}

Finally, we derive the functions
\begin{align*}
G(s)
&=
s^{-1}
(\|\bm{A}^{2}\|
\oplus
\|\bm{h}^{-}\bm{A}\|\|\bm{A}\bm{g}\|
\oplus
\|\bm{h}^{-}\|\|\bm{A}^{2}\bm{g}\|)
\oplus
s^{-2}
\|\bm{A}^{3}\|
\\&=
4s^{-1}
\oplus
5s^{-2},
\\
H(t)
&=
t^{-1}
(\|\bm{A}^{2}\|
\oplus
\|\bm{h}^{-}\bm{A}\|\|\bm{A}\bm{g}\|
\oplus
\|\bm{h}^{-}\|\|\bm{A}^{2}\bm{g}\|)
\oplus
t^{-1/2}
\|\bm{A}^{3}\|^{1/2}
\\&=
4t^{-1}
\oplus
(5/2)t^{-1/2}.
\end{align*}

Observing that the following equalities hold:
\begin{equation*}
\lambda
\oplus
\mu
=
3/2,
\qquad
H(\nu)
=
4\nu^{-1}
\oplus
(5/2)\nu^{-1/2}
=
2,
\end{equation*}
we see that $\lambda\oplus\mu<H(\nu)$. Thus, by Corollary~\ref{C-minxAxxAx-gleqxleqh}, the Pareto frontier for the problem is the set of points $(\alpha,\beta)$ given by the conditions $3/2\leq\alpha\leq2$ and $\beta=4\alpha^{-1}\oplus5\alpha^{-2}$.

Note that for all $\alpha\geq3/2$ the inequalities $4\alpha\geq11/2>5$ hold, which yields $4\alpha^{-1}>5\alpha^{-2}$. As a result, we can reduce the above conditions to those in the form
\begin{equation*}
3/2\leq\alpha\leq2,
\qquad
\beta
=
4\alpha^{-1}.
\end{equation*}

A graphical illustration of the Pareto frontier is given in Fig.~\ref{F-PFSS} (left).

To describe all Pareto-optimal solutions of the problem, we consider the matrix
\begin{equation*}
\alpha^{-1}\bm{A}\oplus\beta^{-1}\bm{1}\bm{1}^{T}\bm{A}
=
\alpha^{-1}\bm{A}\oplus(-4)\alpha\bm{1}\bm{1}^{T}\bm{A}.
\end{equation*}

By applying the condition $3/2\leq\alpha\leq2$, we represent this matrix and derive its square in the form
\begin{equation*}
\left(
\begin{array}{ccc}
1\alpha^{-1} & 2\alpha^{-1} & 2\alpha^{-1}
\\
1\alpha^{-1} & (-2)\alpha & 2\alpha^{-1}
\\
(-3)\alpha & (-2)\alpha & (-2)\alpha
\end{array}
\right),
\qquad
\left(
\begin{array}{ccc}
3\alpha^{-2} & 0 & 4\alpha^{-2}
\\
-1 & 0 & 0
\\
-1 & (-2)\alpha^{2} & 0
\end{array}
\right).
\end{equation*}

Evaluation of the Kleene star matrix yields
%\begin{multline*}
\begin{equation*}
(\alpha^{-1}\bm{A}\oplus\beta^{-1}\bm{1}\bm{1}^{T}\bm{A})^{\ast}
=
%\bm{I}
%\oplus
%\alpha^{-1}\bm{A}\oplus\beta^{-1}\bm{1}\bm{1}^{T}\bm{A}
%\oplus
%(\alpha^{-1}\bm{A}\oplus\beta^{-1}\bm{1}\bm{1}^{T}\bm{A})^{2}
%\\=
\left(
\begin{array}{ccc}
0 & 2\alpha^{-1} & 4\alpha^{-2}
\\
1\alpha^{-1} & 0 & 2\alpha^{-1}
\\
-1 & (-2)\alpha & 0
\end{array}
\right).
%\end{multline*}
\end{equation*}

After calculating the vector
\begin{equation*}
\bm{h}^{-}(\alpha^{-1}\bm{A}\oplus\beta^{-1}\bm{1}\bm{1}^{T}\bm{A})^{\ast}
=
\left(
\begin{array}{rrr}
-1 & 1\alpha^{-1} & 3\alpha^{-2}
\end{array}
\right),
\end{equation*}
we obtain the solution of the problem in the form
\begin{equation*}
\bm{x}
=
\left(
\begin{array}{ccc}
0 & 2\alpha^{-1} & 4\alpha^{-2}
\\
1\alpha^{-1} & 0 & 2\alpha^{-1}
\\
-1 & (-2)\alpha & 0
\end{array}
\right)
\bm{u},
\qquad
\bm{u}^{\prime}
\leq
\bm{u}
\leq
\bm{u}^{\prime\prime},
\qquad
3/2\leq\alpha\leq2,
\end{equation*}
where $\bm{u}=(u_{1},u_{2},u_{3})^{T}$ is the vector of parameters with bounds given by 
\begin{equation*}
\bm{u}^{\prime}
=
\bm{g}
=
\left(
\begin{array}{c}
0
\\
0
\\
0
\end{array}
\right),
\qquad
\bm{u}^{\prime\prime}
=
(\bm{h}^{-}(\alpha^{-1}\bm{A}\oplus\beta^{-1}\bm{1}\bm{1}^{T}\bm{A})^{\ast})^{-}
=
\left(
\begin{array}{c}
1 
\\
(-1)\alpha
\\
(-3)\alpha^{2}
\end{array}
\right).
\end{equation*}

To simplify the solution, we note that the last two columns in the Kleene star matrix are collinear, and hence this matrix can be factored as follows:
\begin{equation*}
\left(
\begin{array}{ccc}
0 & 2\alpha^{-1} & 4\alpha^{-2}
\\
1\alpha^{-1} & 0 & 2\alpha^{-1}
\\
-1 & (-2)\alpha & 0
\end{array}
\right)
=
\left(
\begin{array}{cc}
0 & 2\alpha^{-1}
\\
1\alpha^{-1} & 0
\\
-1 & (-2)\alpha
\end{array}
\right)
\left(
\begin{array}{ccc}
0 & \mathbb{0} & \mathbb{0}
\\
\mathbb{0} & 0 & 2\alpha^{-1}
\end{array}
\right),
\end{equation*}
where we use the notation $\mathbb{0}=-\infty$.

We introduce a new vector of parameters $\bm{v}=(v_{1},v_{2})^{T}$ by the equality
\begin{equation*}
\bm{v}
=
\left(
\begin{array}{ccc}
0 & \mathbb{0} & \mathbb{0}
\\
\mathbb{0} & 0 & 2\alpha^{-1}
\end{array}
\right)
\bm{u}.
\end{equation*}

After turning from $\bm{u}$ to $\bm{v}$, the solution takes the simpler form
\begin{equation*}
\bm{x}
=
\left(
\begin{array}{cc}
0 & 2\alpha^{-1}
\\
1\alpha^{-1} & 0
\\
-1 & (-2)\alpha
\end{array}
\right)
\bm{v},
\qquad
\bm{v}^{\prime}
\leq
\bm{v}
\leq
\bm{v}^{\prime\prime},
\qquad
3/2\leq\alpha\leq2,
\end{equation*}
where the bounds for $\bm{v}$ are derived from $\bm{u}^{\prime}$ and $\bm{u}^{\prime\prime}$ to be
\begin{equation*}
\bm{v}^{\prime}
%=
%\left(
%\begin{array}{ccc}
%0 & \mathbb{0} & \mathbb{0}
%\\
%\mathbb{0} & 0 & 2\alpha^{-1}
%\end{array}
%\right)
%\bm{u}^{\prime}
=
\left(
\begin{array}{c}
0
\\
2\alpha^{-1}
\end{array}
\right),
\qquad
\bm{v}^{\prime\prime}
%=
%\left(
%\begin{array}{ccc}
%0 & \mathbb{0} & \mathbb{0}
%\\
%\mathbb{0} & 0 & 2\alpha^{-1}
%\end{array}
%\right)
%\bm{u}^{\prime\prime}
=
\left(
\begin{array}{c}
1
\\
(-1)\alpha
\end{array}
\right).
\end{equation*}

Moreover, it follows from the two-sided bounds
\begin{equation*}
0
\leq
v_{1}
\leq1,
\qquad
2\alpha^{-1}
\leq
v_{2}
\leq
(-1)\alpha
\end{equation*}
that $1\alpha^{-1}v_{1}\leq2\alpha^{-1}\leq v_{2}$ and $(-1)v_{1}\leq0\leq(-2)\alpha v_{2}$, which allows to write
\begin{equation*}
x_{2}
=
1\alpha^{-1}v_{1}
\oplus
v_{2}
=
v_{2},
\qquad
x_{3}
=
(-1)v_{1}
\oplus
(-2)\alpha v_{2}
=
(-2)\alpha v_{2}.
\end{equation*}

As a result, we represent the Pareto-optimal solution of the problem in the form
\begin{equation}
\bm{x}
=
\left(
\begin{array}{cc}
0 & 2\alpha^{-1}
\\
\mathbb{0} & 0
\\
\mathbb{0} & (-2)\alpha
\end{array}
\right)
\bm{v},
\quad
\left(
\begin{array}{c}
0
\\
2\alpha^{-1}
\end{array}
\right)
\leq
\bm{v}
\leq
\left(
\begin{array}{c}
1
\\
(-1)\alpha
\end{array}
\right),
\quad
3/2\leq\alpha\leq2.
\label{E-x-I-v-I-alpha}
\end{equation}

Note that the vectors $\bm{v}^{\prime}$ and $\bm{v}^{\prime\prime}$ provide the bounds for $\bm{x}$, given by
\begin{equation*}
\bm{x}^{\prime}
=
\left(
\begin{array}{c}
4\alpha^{-2}
\\
2\alpha^{-1}
\\
0
\end{array}
\right),
\qquad
\bm{x}^{\prime\prime}
=
\left(
\begin{array}{c}
1
\\
(-1)\alpha
\\
(-3)\alpha^{2}
\end{array}
\right).
\end{equation*}

We conclude this example by calculating the solutions which correspond to two extreme and one inner points of the Pareto frontier.

First, consider the point $(\alpha,\beta)$ where $\alpha=3/2$ and $\beta=4\alpha^{-1}=5/2$, which corresponds to the best solution with respect to the minimum of the maximum flow time. Then, the solution at \eqref{E-x-I-v-I-alpha} becomes 
\begin{equation*}
\bm{x}
=
\left(
\begin{array}{cc}
0 & 1/2
\\
\mathbb{0} & 0
\\
\mathbb{0} & -1/2
\end{array}
\right)
\bm{v},
\qquad
%\bm{v}^{\prime}
%=
\left(
\begin{array}{c}
0
\\
1/2
\end{array}
\right)
\leq
\bm{v}
\leq
%\bm{v}^{\prime\prime}
%=
\left(
\begin{array}{c}
1
\\
1/2
\end{array}
\right).
\end{equation*}

It is not difficult to verify that in this case we have 
\begin{equation*}
\bm{x}
=
\bm{x}^{\prime}
=
\bm{x}^{\prime\prime}
=
\left(
\begin{array}{cc}
1
\\
1/2
\\
0
\end{array}
\right),
\end{equation*}
and thus the Pareto-optimal solution is the single vector with components
\begin{equation*}
x_{1}=1,
\qquad
x_{2}=1/2,
\qquad
x_{3}=0.
\end{equation*} 

Next, we examine the solution corresponding to $\alpha=2$ and $\beta=4\alpha^{-1}=2$, which is the best with respect to the project makespan. The solution takes the form 
\begin{equation*}
\bm{x}
=
\left(
\begin{array}{cc}
0 & 0
\\
\mathbb{0} & 0
\\
\mathbb{0} & 0
\end{array}
\right)
\bm{v},
\quad
\left(
\begin{array}{c}
0
\\
0
\end{array}
\right)
\leq
\bm{v}
\leq
\left(
\begin{array}{c}
1
\\
1
\end{array}
\right).
\end{equation*}

In scalar form, the solution obtained can be represented as
\begin{equation*}
x_{1}
=
v_{1}\oplus v_{2},
\qquad
x_{2}
=
x_{3}
=
v_{2},
\qquad
0
\leq
v_{1},v_{2}
\leq
1,
\end{equation*}
or, equivalently, using only one parameter $v$, as
\begin{equation*}
v
\leq
x_{1}
\leq
1,
\qquad
x_{2}
=
x_{3}
=
v,
\qquad
0
\leq
v
\leq
1.
\end{equation*}

Finally, consider the inner point of the Pareto frontier with $\alpha=5/3$ and $\beta=4\alpha^{-1}=7/3$. We have the solution
\begin{equation*}
\bm{x}
=
\left(
\begin{array}{cc}
0 & 1/3
\\
\mathbb{0} & 0
\\
\mathbb{0} & -1/3
\end{array}
\right)
\bm{v},
\qquad
\left(
\begin{array}{c}
0
\\
1/3
\end{array}
\right)
\leq
\bm{v}
\leq
\left(
\begin{array}{c}
1
\\
2/3
\end{array}
\right).
\end{equation*}

In scalar form, the solution is written as
\begin{equation*}
x_{1}
=
v_{1}\oplus(1/3)v_{2},
\qquad
x_{2}
=
v_{2},
\qquad
x_{3}
=
(-1/3)v_{2},
\end{equation*}
where $0\leq v_{1}\leq1$ and $1/3\leq v_{2}\leq2/3$.

Another equivalent representation with one parameter $v$ takes the form
\begin{equation*}
(1/3)v
\leq
x_{1}
\leq1,
\qquad
x_{2}
=
v,
\qquad
x_{3}
=
(-1/3)v,
\qquad
1/3\leq v\leq2/3,
\end{equation*}
which can also be rewritten in the usual notation as
\begin{equation*}
v+1/3
\leq
x_{1}
\leq1,
\qquad
x_{2}
=
v,
\qquad
x_{3}
=
v-1/3,
\qquad
1/3\leq v\leq2/3.
\end{equation*}
\end{example}

The optimal solutions obtained are shown in Fig.~\ref{F-PFSS} (right). The thick dot and the big shaded triangle represent the solutions which correspond to the extreme points of the Pareto frontier with $\alpha=3/2$ and $\alpha=2$, whereas the small triangle indicates the solution for the inner point for $\alpha=5/3$.

We now apply Theorem~\ref{T-minxAxxpqx-gxh} to solve problem~\eqref{P-minyixi-yixi-aijxjyi-gixi-yifi}, which can be rewritten in terms of max-plus algebra as
\begin{equation*}
\begin{aligned}
&
\text{minimize}
&&
\left\{
\bigoplus_{1\leq i\leq n}y_{i}x_{i}^{-1},\ 
\bigoplus_{1\leq i\leq n}y_{i}\bigoplus_{1\leq j\leq n}x_{j}^{-1}
\right\};
\\
&
\text{subject to}
&&
\bigoplus_{1\leq j\leq n}a_{ij}x_{j}
=
y_{i},
\\
&&&
g_{i}
\leq
x_{i},
\qquad
y_{i}
\leq
f_{i},
\qquad
i=1,\ldots,n.
\end{aligned}
\end{equation*}

With an additional vector $\bm{y}=(y_{j})$, the problem is represented in vector form as
\begin{equation*}
\begin{aligned}
&
\text{minimize}
&&
\left\{
\bm{x}^{-}\bm{y},\ 
\bm{x}^{-}\bm{1}\bm{1}^{T}\bm{y}
\right\};
\\
&
\text{subject to}
&&
\bm{A}\bm{x}
=
\bm{y},
\quad
\bm{g}
\leq
\bm{x},
\quad
\bm{y}
\leq
\bm{f}.
\end{aligned}
\end{equation*}

We substitute $\bm{y}=\bm{A}\bm{x}$ into the objective functions and the last inequality constraint to eliminate the unknown vector $\bm{y}$. Next, we assume the matrix $\bm{A}$ to be column-regular (to have no zero columns), and the vector $\bm{f}$ to be regular (positive). After application of Lemma~\ref{L-Axleqd} to replace the inequality $\bm{A}\bm{x}\leq\bm{f}$ by the inequality $\bm{x}\leq(\bm{f}^{-}\bm{A})^{-}$, the problem becomes
\begin{equation}
\begin{aligned}
&
\text{minimize}
&&
\{\bm{x}^{-}\bm{A}\bm{x},\ 
\|\bm{x}^{-}\|\|\bm{A}\bm{x}\|\};
\\
&
\text{subject to}
&&
\bm{g}
\leq
\bm{x}
\leq
(\bm{f}^{-}\bm{A})^{-}.
\end{aligned}
\label{P-minxAxxAx-gleqxleqfA}
\end{equation}

To adjust Theorem~\ref{T-minxAxxpqx-gxh} to handle the problem under consideration, we modify the notation \eqref{E-lambda-mu-nu} as
\begin{equation*}
\lambda
=
\bigoplus_{k=1}^{n}
\mathop\mathrm{tr}\nolimits^{1/k}(\bm{A}^{k}),
\qquad
\mu
=
\bigoplus_{k=1}^{n-1}
(\bm{f}^{-}\bm{A}^{k+1}\bm{g})^{1/k},
\qquad
\nu
=
\|\bm{A}\|\oplus\|\bm{f}^{-}\bm{A}\|\|\bm{A}\bm{g}\|,
\end{equation*}
and the notation \eqref{E-G-H} as
\begin{equation*}
\begin{aligned}
G(s)
&=
\bigoplus_{k=1}^{n-1}
s^{-k}
\|\bm{A}^{k+1}\|
\oplus
\bigoplus_{k=1}^{n-2}
s^{-k}
\bigoplus_{\substack{i+j=k\\i,j\geq0}}
\|\bm{f}^{-}\bm{A}^{i+1}\|\|\bm{A}^{j+1}\bm{g}\|,
\\
H(t)
&=
\bigoplus_{k=1}^{n-1}
t^{-1/k}
\|\bm{A}^{k+1}\|^{1/k}
\oplus
\bigoplus_{k=1}^{n-2}
t^{-1/k}
\bigoplus_{\substack{i+j=k\\i,j\geq0}}
\|\bm{f}^{-}\bm{A}^{i+1}\|^{1/k}\|\bm{A}^{j+1}\bm{g}\|^{1/k}.
\end{aligned}
\end{equation*}

We now have the next result.
\begin{corollary}
\label{C-minxAxxAx-gleqxleqfA}
For any column-regular matrix $\bm{A}$, nonzero vector $\bm{g}$ and regular vector $\bm{f}$ such that $\bm{f}^{-}\bm{A}\bm{g}\leq\mathbb{1}$, the following statements hold:
\begin{enumerate}
\item
If $\lambda\oplus\mu\geq H(\nu)$, the Pareto frontier for problem \eqref{P-minxAxxAx-gleqxleqfA} degenerates into the single point $(\alpha,\beta)$ with $\alpha=\lambda\oplus\mu$ and $\beta=\nu$.
\item
Otherwise, the Pareto frontier is the segment given by the conditions
\begin{equation*}
\lambda\oplus\mu
\leq
\alpha
\leq
H(\nu),
\qquad
\beta
=
G(\alpha).
\end{equation*}
\end{enumerate}

All Pareto-optimal solutions are represented in parametric form as
\begin{equation*}
\bm{x}
=
(\alpha^{-1}\bm{A}
\oplus
\beta^{-1}\bm{1}\bm{1}^{T})^{\ast}
\bm{u},
\qquad
\bm{g}
\leq
\bm{u}
\leq
(\bm{f}^{-}\bm{A}
(\alpha^{-1}\bm{A}
\oplus
\beta^{-1}\bm{1}\bm{1}^{T})^{\ast})^{-}.
\end{equation*}
\end{corollary}

\begin{example}
Consider a project with the start-finish relationships, release start and release end time constraints defined by the following matrix and vectors:
\begin{equation*}
\bm{A}
=
\left(
\begin{array}{ccc}
1 & 1 & 2
\\
2 & 1 & \mathbb{0}
\\
\mathbb{0} & 1 & 1
\end{array}
\right),
\qquad
\bm{g}
=
\left(
\begin{array}{c}
0
\\
0
\\
0
\end{array}
\right),
\qquad
\bm{f}
=
\left(
\begin{array}{c}
3
\\
3
\\
2
\end{array}
\right).
\end{equation*}

To apply Corollary~\ref{C-minxAxxAx-gleqxleqfA}, we first perform similar calculations as in the previous example to evaluate matrix powers and vectors involved. Then, we obtain 
\begin{equation*}
\lambda
=
5/3,
\qquad
\mu
=
3/2,
\qquad
\nu
=
2,
\end{equation*}
and dervive the functions
\begin{equation*}
G(s)
=
4s^{-1}
\oplus
5s^{-2},
\qquad
H(t)
=
4t^{-1}
\oplus
(5/2)t^{-1/2}.
\end{equation*}

Further calculation yields $\lambda\oplus\mu=5/3$ and $H(\nu)=4\nu^{-1}\oplus(5/2)\nu^{-1/2}=2$, which show that $\lambda\oplus\mu<H(\nu)$. Therefore, the Pareto frontier is the set of points $(\alpha,\beta)$ given by the conditions $5/3\leq\alpha\leq2$ and $\beta=4\alpha^{-1}\oplus5\alpha^{-2}$.

Since $4\alpha\geq17/3>5$ for all $\alpha\geq5/3$, the inequality $4\alpha^{-1}>5\alpha^{-2}$ is also valid, and thus the above conditions for the Pareto frontier reduce to
\begin{equation*}
5/3\leq\alpha\leq2,
\qquad
\beta
=
4\alpha^{-1}.
\end{equation*}

To describe all Pareto-optimal solutions, we form the Kleene star matrix
\begin{equation*}
(\alpha^{-1}\bm{A}\oplus\beta^{-1}\bm{1}\bm{1}^{T}\bm{A})^{\ast}
=
\left(
\begin{array}{ccc}
0 & 3\alpha^{-2} & 2\alpha^{-1}
\\
2\alpha^{-1} & 0 & 4\alpha^{-2}
\\
(-2)\alpha & -1 & 0
\end{array}
\right).
\end{equation*}

After decomposition of the matrix in the form
\begin{equation*}
\left(
\begin{array}{ccc}
0 & 3\alpha^{-2} & 2\alpha^{-1}
\\
2\alpha^{-1} & 0 & 4\alpha^{-2}
\\
(-2)\alpha & -1 & 0
\end{array}
\right)
=
\left(
\begin{array}{cc}
0 & 3\alpha^{-2}
\\
2\alpha^{-1} & 0
\\
(-2)\alpha & -1
\end{array}
\right)
\left(
\begin{array}{ccc}
0 & \mathbb{0} & 2\alpha^{-1}
\\
\mathbb{0} & 0 & \mathbb{0}
\end{array}
\right),
\end{equation*}
we finally represent the solution of the problem as
\begin{equation*}
\bm{x}
=
\left(
\begin{array}{cc}
0 & 3\alpha^{-2}
\\
2\alpha^{-1} & 0
\\
(-2)\alpha & -1
\end{array}
\right)
\bm{v},
\qquad
5/3\leq\alpha\leq2,
\end{equation*}
where the vector of parameters satisfies the condition
\begin{equation*}
\left(
\begin{array}{c}
2\alpha^{-1}
\\
0
\end{array}
\right)
\leq
\bm{v}
\leq
\left(
\begin{array}{c}
(-1)\alpha
\\
1
\end{array}
\right).
\end{equation*}

It follows from the bounds for the vector $\bm{v}$ that $v_{1}\geq2\alpha^{-1}\geq3\alpha^{-2}\geq3\alpha^{-2}v_{2}$ and $(-2\alpha)v_{1}\geq0\geq(-1)v_{2}$, which allows to write
\begin{equation*}
x_{1}
=
v_{1}
\oplus
3\alpha^{-2}v_{2}
=
v_{1},
\qquad
x_{3}
=
(-2)\alpha v_{1}
\oplus
(-1)v_{2}
=
(-2)\alpha v_{1}.
\end{equation*}

As a result, the Pareto-optimal solution becomes
\begin{equation*}
\bm{x}
=
\left(
\begin{array}{cc}
0 & \mathbb{0}
\\
2\alpha^{-1} & 0
\\
(-2)\alpha & \mathbb{0}
\end{array}
\right)
\bm{v},
\end{equation*}
where $\bm{v}$ and $\alpha$ are given by the same conditions as before.

Using scalar representation yields the solution in the form
\begin{equation*}
x_{1}
=
v_{1},
\qquad
x_{2}
=
2\alpha^{-1}v_{1}
\oplus
v_{2},
\qquad
x_{3}
=
(-2)\alpha v_{1}.
\end{equation*}

Finally, turning to the usual notation, we rewrite the solution as
\begin{equation*}
x_{1}
=
v_{1},
\qquad
x_{2}
=
\max\{v_{1}-\alpha+2,v_{2}\},
\qquad
x_{3}
=
v_{1}+\alpha-2,
\end{equation*}
where $2-\alpha\leq v_{1}\leq\alpha-1$, $0\leq v_{2}\leq1$ and $5/3\leq\alpha\leq2$.
\end{example}

\section{Conclusions}

Bi-objective project scheduling problems without cost dependencies and resource requirements have been considered. Given temporal constraints that include start-finish precedence relations, release start times, release end times and deadlines, the problems are to develop a schedule that minimizes both the maximum flow-time over all activities and the project makespan. The solution of such problems can be used as an auxiliary tool in solving more general project scheduling problems, and is of independent interest.

We have represented the problems as a bi-objective optimization problem in terms of tropical (idempotent) mathematics which concerns with the theory and applications of idempotent semirings and semifields. By using methods of tropical optimization, we have obtained complete Pareto-optimal solutions of the problems in an exact analytical form suitable for formal analysis and computations with polynomial time.

The results obtained demonstrate that application of tropical optimization techniques allows one to obtain analytic solutions to real-world problems that, in practice, are solved numerically by using various computational procedures and do not have exact explicit solutions available. The proposed solutions may efficiently serve to complement and supplement the existing numerical approaches, and become the only solution when the algorithmic solutions are infeasible or impossible to implement. 

Further investigation can include the derivation of solutions to bi-objective and multiobjective problems with new criteria, including minimum deviation of the finish times of activities from due dates, minimum deviation of start or finish times of activities, as well as with additional constraint, such as due-dates, start-start and finish-start precedence relationships.

%\section*{Acknowledgments}
%This work was supported in part by the Russian Foundation for Basic Research (grant number 18-010-00723).

\bibliographystyle{abbrvurl}

\bibliography{Tropical_optimization_technique_in_bi-objective_project_scheduling_under_temporal_constraints}

\end{document}